\documentclass{amsart}
\usepackage{amssymb}
\usepackage{thmtools}
\usepackage{relsize}
\usepackage[T1]{fontenc}

\usepackage{constants} 

\usepackage{tikz}
\usepackage{tikz-cd}

\usepackage{hyperref}
\usepackage[capitalize]{cleveref}
\usepackage{autonum}

\newtheorem{theorem}{Theorem}
\newtheorem*{theorem*}{Theorem}
\newtheorem{corollary}[theorem]{Corollary}
\newtheorem{hypothesis}[theorem]{Hypothesis}
\newtheorem{proposition}[theorem]{Proposition}
\newtheorem{lemma}[theorem]{Lemma}

\theoremstyle{definition}
\newtheorem{definition}[theorem]{Definition}
\crefname{definition}{Definition}{Definitions}

\theoremstyle{remark}
\newtheorem{remark}[theorem]{Remark}

\newconstantfamily{c}{symbol=c}

\newcommand{\OK}{\mathcal{O}_K}

\newcommand{\prim}{prim}
\newcommand{\KR}{K_\mathbb{R}}

\newcommand{\Q}{\mathbb{Q}}
\newcommand{\diff}{\,\mathrm{d}}
\newcommand{\hecke}[2]{ \underset{( #1, #2 )}{\leftarrow} }

\DeclareMathOperator{\SL}{SL}
\DeclareMathOperator{\M}{M}

\DeclareMathOperator{\GL}{GL}

\DeclareMathOperator{\argmin}{argmin}
\DeclareMathOperator{\vol}{vol}

\DeclareMathOperator{\card}{\texttt{\#}}

\DeclareMathOperator{\rank}{rk}
\DeclareMathOperator{\spant}{span}
\DeclareMathOperator{\ind}{\mathlarger{\mathbf{1}}}
\DeclareMathOperator{\supp}{supp}
\DeclareMathOperator{\N}{N}
\DeclareMathOperator{\Tr}{tr}

\DeclareMathOperator{\Gr}{\mathbf{Gr}}

\begin{document}
\title{Integral matrices of fixed rank over number fields}
\author[N. Gargava]{Nihar Gargava}
\author[V. Serban]{Vlad Serban}
\author[M. Viazovska]{Maryna Viazovska}
\author[I. Viglino]{Ilaria Viglino}

\address{N. Gargava, Orsay Institute of Mathematics, Paris-Saclay University, France}
\email{nihar.gargava@universite-paris-saclay.fr}
\address{V. Serban, Department of Mathematics, New College of Florida, U.S.A.}
\email{vserban@ncf.edu}
\address{M. Viazovska, Section of Mathematics, EPFL, Switzerland}
\email{maryna.viazovska@epfl.ch}
\address{I. Viglino, Section of Mathematics, EPFL, Switzerland}
\email{ilaria.viglino@epfl.ch}

\begin{abstract}
    We prove an asymptotic formula for the number of fixed rank matrices with integer coefficients over a number field $K/\mathbb{Q}$ and bounded norm. As an application, we derive an approximate Rogers integral formula for discrete sets of module lattices obtained from lifts of algebraic codes. This in turn implies that the moment estimates of \cite{GSV23},{ which inform the behavior of short vectors in sets of random lattices, }also carry through for large enough discrete sets of module lattices. 
\end{abstract}
\maketitle

\section{Introduction}
\label{se:intro}

We start by revisiting a fundamental counting result for integral matrices of fixed rank by Y. Katznelson \cite{K1994}. Fix integers $n > m \geq k \geq 1$. The main result of \cite{K1994} proves the following asymptotic counts: 
\begin{theorem}
    Let $f: \M_{n \times m}( \mathbb{R}) \rightarrow  \mathbb{R}$ be the indicator function of an origin-centered unit ball in the $l^{2}$-norm $\| \cdot \| : \M_{n \times m}(\mathbb{R})\cong \mathbb{R}^{nm} \rightarrow   \mathbb{R}$. 
    Then, for some constants $\Cl[c]{main}, \Cl[c]{error} > 0$ that depend on $n,m,k$ but not on $T \geq 1$, one has 
    \begin{equation}
        \sum_{\substack{A \in \M_{ n \times m }(\mathbb{Z}) \\ \rank(A) = k }} f( \tfrac{1}{T} A)  =  \Cr{main} \cdot T^{ k n} \cdot (1 + \varepsilon),
        \label{eq:main_eq1}
    \end{equation}
    where 
    \begin{equation}
        |\varepsilon| \leq \Cr{error} \cdot T^{-1}.
        \label{eq:error_term1}
    \end{equation}
    \label{th:main1}
\end{theorem}
This theorem solves an interesting counting problem. Indeed, because $f$ is the indicator function of a unit $l^{2}$-ball we have: 
\begin{equation}
    \sum_{\substack{A \in \M_{ n \times m }(\mathbb{Z}) \\ \rank(A) = k }} f( \tfrac{1}{T} A)  
    = \card \{A \in  \M_{n \times m}(\mathbb{Z}) \mid  \rank(A) = k, \|A\| \leq T \}.
\end{equation}
This result of Katznelson has motivated many subsequent refinements and generalizations, see for example \cite{MOS24,MOS25,B16}.

The first main result of our paper establishes a natural number-theoretic generalization of Katznelson's counting result. 
Let $K$ be a number field of degree $\deg K = d$,
and let $\OK$ denote the ring of integers of the number field. 
For notational simplicity, we will abbreviate $\KR = K \otimes \mathbb{R}$. 
We consider the analogous counting problem over number fields, so that we work with matrices with $\OK$-entries. Moreover, in the spirit of similar results in the geometry of numbers, it seems natural to allow summation over more general functions such as compactly supported smooth functions or indicator functions of bounded convex sets on the space of matrices $\M_{n \times m}(\KR)$. For a class of so-called admissible functions which includes the examples above, we show:
\begin{theorem}
    Let $f: \M_{n \times m}( \KR) \rightarrow  \mathbb{R}$ be an admissible function (see \cref{hy:admissible}). 
    Then, for some constants $\Cr{main} \in \mathbb{R}$ and $\Cr{error}>0$ that depend on $K,n,m,k,f$ but not on $T \geq 2$, one has 
    \begin{equation}
        T^{-knd}\sum_{\substack{A \in \M_{ n \times m }(\OK) \\ \rank(A) = k }} f( \tfrac{1}{T} A)  =  \Cr{main} + \varepsilon,
        \label{eq:main_eq}
    \end{equation}
    where 
    \begin{equation}
        |\varepsilon| \leq \Cr{error} \cdot T^{-1} \log T.
        \label{eq:error_term}
    \end{equation}
    \label{th:main}
    Moreover, unless $d=k=1$ and $m=n-1$, the $\log T$ in~\cref{eq:error_term} can be dropped (cf.~\cref{ss:log_term} about the $\log T$ factor).
\end{theorem}

\subsection{Connection to Rogers's integral formula}
It seems furthermore reasonable to expect that the leading constant $\Cr{main} \in \mathbb{R}$ in the main factor of the asymptotic formula in \cref{eq:main_eq} carries some arithmetic-geometric meaning. To that end, we highlight a striking connection of the leading constant in both Katznelson's and our more general work to Rogers's integration formula in the geometry of numbers \cite{GSV23,Rogers55, Rogers58}
. 

Let us explain. Observe that for any $A \in \M_{n \times m }(\OK)$ 
such that $\rank(A) =k$, one can perform the rank factorization and write 
\begin{equation}
    A = C \cdot D,
\end{equation}
where $C \in \M_{n \times k}(K)$, $D \in \M_{k \times m}(K)$ and $\rank(C) = \rank(D) = k$. Assuming that $D$ is in row-reduced echelon form of maximal rank, the choice of $C$ and $D$ is unique. 
For brevity, for the rest of the paper we will abbreviate ``echelon'' for  ``row-reduced echelon of maximal rank''.
We then observe the following decomposition, which enables the connection to Rogers's formula:
\begin{align}
 & \{A \in \M_{n \times m }(\OK) \mid \rank(A) = k \}   \\
    = &  \bigsqcup_{ \substack{ D \in \M_{k \times m}(K) \\ D \text{ echelon }}}
    \{C  \in \M_{n \times k }(K) \mid  C \cdot D \in \M_{ n \times m }(\OK), \rank(C) =k \} \cdot D.
    \label{eq:bijection}
\end{align}

For an echelon matrix $D \in \M_{k \times m }(K)$ the denominator $\mathfrak{D}(D) \in \mathbb{Z}_{\geq 1}$ is given by the following index:
\begin{equation}
    \mathfrak{D}(D) = [ \OK^{k} :\{v \in  \OK^{k} \mid D^{T} v \in \OK^{m} \}].
    \label{de:denonimator}
\end{equation}
We show in this paper that the constant $\Cr{main}$ in \cref{eq:main_eq} is precisely given by
\begin{equation}
    \Cr{main} = \sum_{ \substack{ D \in \M_{k \times m}(K) \\ D \text{ echelon}} } \mathfrak{D}(D)^{-n}\int_{x \in \M_{n \times k}(\KR)} f(x D )  \diff x,
    \label{eq:value_of_c}
\end{equation}
where the integral is over the Euclidean structure on $\M_{n \times k}(\KR)$ given by \cref{eq:norm}. 
Convergence of the infinite sum on the right-hand side when $n > m \geq k$ in fact follows from Schmidt's seminal work \cite{S1967} on rational points of Grassmannian varieties over number fields, as explained in \cref{se:schmidt}.

Heuristically, \cref{eq:value_of_c} can be understood as follows: using the bijection in \cref{eq:bijection}, 
we can at least formally write 
\begin{equation}
    T^{-knd}\sum_{\substack{A \in \M_{ n \times m }(\OK) \\ \rank(A) = k }} f( \tfrac{1}{T} A)  
    = \sum_{\substack{D \in \M_{k \times m}(\KR) \\ D \text{ echelon} }} \Big[ T^{-knd}\sum_{\substack{C \in \M_{n \times k}(K)  \\ C \cdot D   \in \M_{n \times m} (\OK) \\ \rank(C) = k }} f(\tfrac{1}{T} C\cdot D ) \Big].
    \label{eq:sums_distributes}
\end{equation}
Now it can be easily argued that the contribution of each echelon matrix $D$ approximates a Riemann integral. More precisely, as $T \rightarrow  \infty$ one can observe that the Euclidean measure on $\M_{n \times k}(\KR)$ is chosen conveniently so that
\begin{equation}
    \sum_{\substack{C \in \M_{n \times k}(K)  \\ C \cdot D   \in \M_{n \times m} (\OK) \\ \rank(C) = k }} 
    {T^{-knd}}
    f(\tfrac{1}{T} C\cdot D )  \rightarrow
    \mathfrak{D}(D)^{-n}
    \int_{\M_{n \times k}(\KR) } f(xD) \diff x.
\end{equation}

Therefore, one knows that \cref{eq:value_of_c} yields the obvious candidate value for $\Cr{main}$. 
However, it is nontrivial to show that the exchange of limits is possible in \cref{eq:sums_distributes}. 
For example, each of the inner terms might have an error which accumulates over infinitely many $D$ in the outer sum. 
Proving our main theorem essentially amounts to showing that such issues do not arise.

\par 
{
    Interestingly, Katznelson \cite{K1994} uses some formulas of Terras \cite{T12}
    to express the constant $c_{1}$ in terms of special values of Koecher zeta functions. Combining this with our observations leads to interesting equalities. For instance, in the case of $K=\mathbb{Q}$ and $f$ the indicator function of the origin-centered unit ball in $\mathbb{R}^{nm}$, one obtains the formula: 

    \begin{equation}
        \sum_{ \substack{ D \in \M_{k \times m}(K) \\ D \text{ echelon}} } \mathfrak{D}(D)^{-n}\int_{x \in \M_{n \times k}(\KR)} f(x D )  \diff x=\frac{V(nk)\cdot Z_{k,m-k}(I, n/2)}{\zeta(n)\cdot \zeta(n-1)\cdots \zeta(n-k+1)}
        \label{eq:correct_usage}
    \end{equation}
    where we denote henceforth by $V(s)$ the unit ball's volume in $\mathbb{R}^s$ and where for a positive symmetric matrix $X\in M_{m\times m}(\mathbb{R})$ with $\Re(s)>m/2$ the Koecher zeta function is defined by 
$$Z_{k,m-k}(X, s)=\sum_{L\in \mathbb{Z}^{m\times k}/\GL_k(\mathbb{Z})}\det(L^tXL)^{-s}. $$}

\subsection{Motivation from coding theory and cryptography}

Consider an $\OK$-module $\Lambda \subseteq K^{n} \otimes \mathbb{R}$ of $\OK$-rank $n$.
For a prime ideal $\mathcal{P} \subseteq \OK$ and $1 \leq r \leq n$, we will now define a $(\mathcal{P},r)$-Hecke neighbor $\Lambda'$ to be a lattice given by the following construction. 

\begin{definition}
    \label{de:heckenhbr}
    Let $k_{\mathcal{P}} = \OK/\mathcal{P}$ be the residue field 
    of $\mathcal{P}$ and let $\N(\mathcal{P}) = \card k_{\mathcal{P}}$ be the ideal norm of $\mathcal{P}$.
    Let $\pi_{\mathcal{P}}$ be the ``modulo $\mathcal{P}$'' reduction map given as
    \begin{equation}
        \pi_{\mathcal{P}}:\Lambda \rightarrow \Lambda/\mathcal{P}\Lambda \simeq  k_{\mathcal{P}}^{n}.
        \label{eq:modp}
    \end{equation}
    One then says that a lattice $\Lambda' \subseteq \KR^{n}$ is a $(\mathcal{P},r)$-Hecke neighbor of $\Lambda$ if
    for some $r$-dimensional $k_{\mathcal{P}}$-subspace $V \subseteq k_{\mathcal{P}}^{n}$,
    \begin{equation}
        \Lambda'= \N(\mathcal{P})^{- \tfrac{1}{d}\left(1-\tfrac{r}{n}\right)} \pi_{\mathcal{P}}^{-1}(V).
        \label{eq:defi_of_hecke_nbr}
    \end{equation}
    We will abbreviate $ \Lambda' \hecke{\mathcal{P}}{r} \Lambda$ to say that $\Lambda'$ is a $(\mathcal{P},r)$-Hecke neighbor of $\Lambda$.
\end{definition}

The scaling factor in front of $\pi_{\mathcal{P}}^{-1}(V)$ in \cref{eq:defi_of_hecke_nbr}
ensures that $\vol( \KR^{n}/\Lambda) = \vol( \KR^{n}/\Lambda')$.

Given a lattice $\Lambda \subseteq \KR^{n}$, the number of lattices $\Lambda' \subseteq \KR^{n}$ that are $(\mathcal{P},r)$-Hecke neighbors of $\Lambda$ is exactly the cardinality of the Grassmannian variety $\Gr(r,k_{\mathcal{P}}^{n})$ over the finite field $k_{\mathcal{P}}$.
Due in part to this finiteness property, such constructions of lattices have drawn considerable interest in algorithmic applications of lattices.
They are referred to as ``lifts of codes'' \cite{GS21}
or ``Construction A'' lattices in the coding theory literature \cite{CS,VKO14}; 
see also the literature on $q$-ary lattices, for example \cite{C18}.
In lattice-based cryptography, such Hecke neighbors appear in ``worst-case to average-case'' reductions \cite{AD15}.

Rogers \cite{Rog47} was perhaps the first to study random $(p,r)$-Hecke neighbors
for the case of $\mathbb{Z}^{n} \subseteq \mathbb{R}^{n}$, an integer prime $p$ and $1 < r < n$. 
His key observation is that $(p,r)$-Hecke neighbors satisfy Siegel's mean value theorem on average as $p \rightarrow \infty$. That is, for any admissible test function $f:\mathbb{R}^{n} \rightarrow  \mathbb{R}$, one has the convergence of expected values for lattice sums
\begin{equation}
    \mathbb{E}_{\Lambda \hecke{p}{r} \mathbb{Z}^{n} }\Big( \sum_{v \in \Lambda \setminus \{0\}} f(v) \Big) 
    \rightarrow \int_{\mathbb{R}^{n}} f(x) \diff x,
    \label{eq:rogers_conv}
\end{equation}
for $p \rightarrow  \infty$. 
The modern way to understand this convergence is via Siegel's mean value theorem. Writing $\mu_{\Pr}$ for the Haar-probability measure on $\SL_{n}(\mathbb{R})/\SL_{n}(\mathbb{Z})$, it states that
\begin{equation}
    \int_{\SL_{n}(\mathbb{R})/\SL_{n}(\mathbb{Z})} \Big( \sum_{v \in \Lambda \setminus \{0\}} f(v) \Big) 
    \diff \mu_{\Pr}
    = \int_{\mathbb{R}^{n}} f(x) \diff x.
\end{equation}
Moreover, due to the more modern work \cite{EOL01}, the fact that $(p,r)$-Hecke neighbors of a fixed lattice (say $\mathbb{Z}^{n}$) equidistribute in the space $\SL_{n}(\mathbb{R})/\SL_{n}(\mathbb{Z})$ as $p\to \infty$ is well-established. Therefore \cref{eq:rogers_conv} should be expected to hold \footnote{However, proving this from the equidistribution of $(p,r)$-Hecke neighbors would probably require a stronger version of \cite[Theorem 1.6]{EOL01}. We would need a pointwise convergence statement for Hecke operators for non-compactly supported functions (since $\Lambda \mapsto \sum_{v \in\Lambda \setminus \{0\}} f(v)$ is not compactly supported). We would like to thank an anonymous referee for pointing this out.}.

However, a direct proof of \cref{eq:rogers_conv} as explained in \cite{Rog47} is elementary. All one has to do is rewrite 
\begin{align}
& \mathbb{E}_{\Lambda \hecke{p}{r} \mathbb{Z}^{n} }\Big( \sum_{v \in \Lambda \setminus \{0\}} f(v) \Big)  \\
= 
& \Big(\sum_{v \in p \mathbb{Z}^{n} \setminus \{0\}} f(p^{-\left(1-\frac{r}{n}\right)} v)\Big) + 
\mathbb{E}_{V \in \Gr(r,k_{p}^{n})} 
\Big(\sum_{v \in \mathbb{Z}^{n}  \setminus p \mathbb{Z}^{n}} \ind_{\pi_{p}(v) \in V} f ( p^{- \left( 1-\frac{r}{n}  \right) v })\Big).
\end{align}

One then observes that the first term must converge to $0$. Indeed, since $p \cdot p^{-(1-r/n)} = p^{r/n} \rightarrow \infty$, as $p$ grows the first sum will not contain any points in the support of $f$. On the other hand, after substituting the cardinality $\card \Gr(r,k_{p}^{n})$, one can show that the second term converges towards a Riemann integral approximating $\int_{\mathbb{R}^{n}}f(x) \diff x$ as $p \rightarrow  \infty$.

For $2\leq m < n$, one may hope to generalize this elementary proof to $m$-th moments by considering the expectation
\begin{equation}
    \label{eq:higher_moments_hecke}
    \mathbb{E}_{\Lambda \hecke{p}{r}  {\mathbb{Z}^{n}}} \Big[ \Big( \sum_{v \in \Lambda \setminus \{0\}}  f(v) \Big)^{m}\Big].
\end{equation}
However, one quickly runs into having to prove \cref{th:main} with $K = \mathbb{Q}$. 
We therefore show in \cref{se:lifts_of_codes} how our results allow us to evaluate over arbitrary number fields $K$  moments as in \cref{eq:higher_moments_hecke} as the norm of the prime $\mathcal{P}$ goes to infinity. { In particular, we obtain as a consequence of Theorem \ref{th:higher_moments}: 
\begin{theorem}\label{thm:mainapplication}
Let $n \ge 2$, $m \in \{1,\dots,n-1\}$ and $r$ be chosen as any number in $\{m,m+1,\dots,n-1\}$ satisfying $1-\frac{r}{n} < \frac{1}{m}$
(e.g. $r=n-1$ works always)
. Let $f:K_{\mathbb{R}}^{n} \rightarrow \mathbb{R}$ be a function satisfying \cref{hy:admissible}. 
As $\N(\mathcal{P}) \rightarrow \infty$ we have the convergence 
{\small \begin{equation}
\mathbb{E}_{\Lambda \hecke{\mathcal{P}}{r}  {{\OK}^{n}}} \Big[ \Big( \sum_{v \in \Lambda \setminus \{0\}}  f(v) \Big)^{m}\Big]\to \int_{\SL_n(\KR)/\SL_n(\OK)}\Big( \sum_{v \in \Lambda \setminus \{0\}}  f(v) \Big)^{m}\diff \mu_{\Pr}.
\end{equation}}
In other words, moments over the discrete spaces of Hecke neighbors approximate moments for the full space of Haar-random free $\OK$-modules of unit covolume. 
\end{theorem}
As an immediate corollary, we therefore deduce that the moment estimates of \cite{GSV23}, which control the behavior of short vectors in Haar-random number field lattices, also apply for primes $\mathcal{P}$ of large enough norm to the discretized sets of Hecke neighbors of an $\OK$-lattice. }

\subsection{The case of \texorpdfstring{$k=1,d=1,n=m+1$}{k=1,d=1,n=m+1}}
\label{ss:log_term}
There appears to be a technical gap in the proof of \cite{K1994} for the case of $k=1,d=1,n=m+1$ that was perhaps overlooked by the author but can be fixed as follows.  

Let $f:\M_{n \times (n-1)}(\mathbb{R}) \rightarrow  \mathbb{R}$ be the indicator function of a unit ball.
In this case, one wants to sum for $T \geq 1$
\begin{align}
    \sum_{\substack{A \in \M_{n \times (n-1)}(\mathbb{Z}) \\ \rank A = 1}} f( \tfrac{1}{T}A)
    =
    \frac{1}{2}\sum_{\substack{v \in \M_{1 \times (n-1)}(\mathbb{Z}) \setminus \{0\} \\ \gcd(v_{1},\dots,v_{n-1}) = 1} } 
    \card \{ w \in \M_{n \times 1}(\mathbb{Z}) \setminus \{0\} \mid \|w v\| \leq T \} 
\end{align}
The factor $1/2$ comes from the symmetry of $w \cdot v = (-w) \cdot (-v)$.
Now it turns out for a column matrix $w$ and a row matrix $v$ one has $\|wv\|=\|w\|\|v\|$. Hence, the sum becomes
\begin{align}
    \frac{1}{2}\sum_{\substack{v \in \M_{1 \times (n-1)}(\mathbb{Z}) \setminus \{0\} \\ 1 \leq \|v\| \leq T , \gcd(v)=1} } 
    \card \{ w \in \M_{n \times 1}(\mathbb{Z}) \setminus \{0\} \mid \|w\| \leq T \|v\|^{-1} \} 
\end{align}
If we bound for some constant $\Cl[c]{gauss1}>0$ the cardinality of the set 
\begin{equation}
    \card \{ w \in \M_{n \times 1}(\mathbb{Z}) \setminus \{0\} \mid \|w\| \leq T \|v\|^{-1} \} \leq V(n)\frac{T^{n}}{\|v\|^{n}} + \Cr{gauss1}{ \frac{T^{n-1}}{\|v\|^{n-1}}},
\end{equation}
then we observe that $\sum_{v \in \M_{1 \times (n-1)}(\mathbb{Z}), \gcd(v)=1} \|v\|^{-n}$ is finite whereas 
the second term must contribute to a $\sim \log T$ factor. However, if we use any non-trivial bound on the Gauss circle problem in $n \geq 2$ dimensions, we get
\begin{equation}
    \card \{ w \in \M_{n \times 1}(\mathbb{Z}) \setminus \{0\} \mid \|w\| \leq T \|v\|^{-1} \} \leq V(n)\frac{T^{n}}{\|v\|^{n}} + \Cl[c]{gauss2}{ \frac{T^{n-1- \delta} }{\|v\|^{n-1 - \delta}}} ,
\end{equation}
for some $\delta>0$ and no $\log T$ term is introduced in the error term. This is because by summation by parts, one gets
\begin{equation}
    \sum_{\substack{v \in \M_{1 \times (n-1) (\mathbb{Z})} \\ \gcd(v)=1, \|v\|\leq T} } \|v\|^{-(n-1-\delta)}  \leq \Cl[c]{gauss3} T^{\delta}.
\end{equation}

In the setting of~\cref{th:main}, we do not restrict $f$ to be the indicator function of some $l^{2}$-ball and therefore bounds from the Gauss circle problem do not necessarily apply. Hence, the $\log T$ term cannot be removed for the case of $d=1,k=1,n=m+1$ unless we change our~\cref{hy:admissible}.

\subsection{Some comments on the proof}

Our proof is heavily inspired from the combinatorial arguments of Katznelson \cite{K1994}. 
In particular, we use an analogous Minkowski reduction theory extended to number fields. Independently, Seungki Kim established in \cite[Prop 3.2]{K19} a result similar to Theorem \ref{thm:mainapplication} for rank $k=1$ in the ad\`elic formulation and with some restrictions on the primes $\mathcal{P}$ allowed.

In our work, the focus in \cref{th:main} is on asymptotics for $T \rightarrow  \infty$.
We do not explore the variation of $\Cr{error}$ in terms of the number field $K,f$ and the integer constants $k,m,n$.
However, a so inclined reader should be able to chase through our constants $c_{1},c_{2},\dots$ to understand the dependence on these parameters.
We have made little effort to optimize this dependence, and it seems unlikely that the best route to optimize $\Cr{error}$ is through this combinatorial approach.

\section*{Acknowledgements}

We would like to thank Andreas Str\"ombergsson for comments on our preprint \cite{GSV23} that motivated the writing of much of this article. The writing has also benefited from some discussions with Phong Nguyen, Thomas Espitau and Seungki Kim.  
We owe gratitude to Fabian Gundlach for pointing out an issue with~\cref{hy:admissible} and also suggesting a speedy resolution. Finally, we thank the anonymous referees for their careful reading, for spotting a number of omissions and for numerous contributions to improve the quality of this text.

This research funded in part by the SNSF Project funding 184927 titled "Optimal configurations in multidimensional spaces". 
NG acknowledges funding from the ERC Grant
101096550 titled ``Integrating Spectral and Geometric data on Moduli Space'' and 
from the SNSF grant 225437 titled "Random Geometry with Arithmetic Constraints".
MV acknowledges funding from the SNSF grant 215337 titled ``Sphere packing, Energy minimization and Fourier interpolation''.

AI tools were used to get feedback on the writing, but the paper is entirely human written.
However, with the recent rise of AI-assisted mathematics, we were inspired to create an autoformalization
of this project. The main results~\cref{th:main} and~\cref{thm:mainapplication} are now machine verified by the Palomar registry; see the accompanying Lean formalization~\cite{IntegralMatrices}.
The formalization assumes Schmidt's theorem (\cref{th:schmidt}), which we did not formalize.

%
%
%

\section{Preliminaries and Notations}

\subsection{Lattices}

\begin{definition}
    We define a lattice in a Euclidean space $V$ to be a discrete $\mathbb{Z}$-submodule $\Lambda \subseteq V$. 
    A lattice has finite covolume in $V$ if 
    $\vol(V/\Lambda) < \infty$.
\end{definition}
\begin{remark}
    Given any ambient Euclidean space, a lattice $\Lambda \subseteq V$ has finite covolume in $\Lambda \otimes_{\mathbb{Z}} \mathbb{R}$.
\end{remark}

\begin{definition}
    \label{de:height_definition}
    For any Euclidean space $V$ and for any discrete $\mathbb{Z}$-submodule $\Lambda \subseteq V$, we define the height 
    of $\Lambda$ by
    \begin{equation}
        H(\Lambda)  = \vol(\Lambda \otimes \mathbb{R} /\Lambda)
    \end{equation}
    taken with respect to the restriction of the norm to $ \Lambda \otimes \mathbb{R} \subseteq V$. 
\end{definition}

\begin{definition}
    \label{de:primitive}
    Let $\Lambda \subseteq V$ be a lattice and let $\Lambda' \subseteq \Lambda$ be a sublattice.
    We call the lattice $\Lambda$-primitive if $(\Lambda' \otimes \mathbb{Q}) \cap \Lambda = \Lambda'$.

    Most of the time, we will skip mentioning $\Lambda$ when the context is clear.
\end{definition}

\subsection{Covering radius and Voronoi domain}

\begin{definition}
    For a lattice $\Lambda \subseteq V$ in Euclidean space $V$, we denote by $\rho(\Lambda)$ the covering radius of $\Lambda$ defined as 
    \begin{equation}
        \rho(\Lambda) = \max_{x \in {\Lambda \otimes \mathbb{R}}} \min _{v \in \Lambda} \|x-v\|.
    \end{equation}
\end{definition}

\begin{definition}
    Given a lattice $\Lambda \subseteq V$ in a Euclidean space $V$, one defines a Voronoi domain $F \subseteq \Lambda \otimes \mathbb{R}$ as
    \begin{equation}
        F = \{ x \in \Lambda \otimes \mathbb{R} \mid \|x\| \leq \|x + v\| \text{ for all } v \in \Lambda\}.
    \end{equation}
\end{definition}
One has the following properties of the Voronoi domain.
\begin{lemma}
    \label{le:voronoi}
    \begin{enumerate}

        \item  We have $F + \Lambda = \Lambda \otimes \mathbb{R}$,
        \item We have $\vol(F) = H(\Lambda)$,
        \item  One has $F \subseteq B_{0}(\rho(\Lambda))$, where $B_{0}(\rho(\Lambda))$ is the ball of radius $\rho(\Lambda)$ and center at $0$.
    \end{enumerate}
    \begin{proof}
        Standard facts.
    \end{proof}

\end{lemma}

\begin{lemma}
    \label{le:ballvol}
    For a lattice $\Lambda$ in a Euclidean space $V$, one has that for any $T > 0$
    \begin{equation}
        \card\{v \in \Lambda \mid \|v\| \leq T \} \leq \Cl[c]{ballvol} (T + \rho(\Lambda))^{r} H(\Lambda)^{-1},
    \end{equation}
    where $r= \rank_{\mathbb{Z}} \Lambda$ and $\Cr{ballvol}$ is a constant depending only on $r$.
\end{lemma}
\begin{proof}
    This is a volume argument. We take a Voronoi domain $F \subseteq \Lambda \otimes \mathbb{R}$. Then the set 
    \begin{equation}
        F + \{v \in \Lambda \mid \|v\| \leq  T\} \subseteq\{v \in  \Lambda \otimes \mathbb{R} \mid \|v\| \leq T + \rho(\Lambda) \}.
    \end{equation}
\end{proof}

For the purpose of this article, we will assume that $\Cr{ballvol}>0$ is large enough to work for all $r \leq dn^{3}$ necessary for our purposes.

\subsection{Minkowski and Hadamard}

The following is an important lemma due to Minkowski. 
\begin{lemma}
    \label{le:minkowski1}
    Let $\Lambda \subseteq V$ be a lattice in a Euclidean space $V$ whose $\mathbb{Z}$-rank is $r$. 
    Then, for some vector $v \in \Lambda \setminus \{0\}$, we have 
    \begin{equation}
        \|v\|\leq  \Cl[c]{hermite} H(\Lambda)^{\tfrac{1}{r}}  
    \end{equation}
    for some $\Cr{hermite}>0$ depending on $r$.
\end{lemma}
Although the constant $\Cr{hermite}$ depends on the $\mathbb{Z}$-rank $r$ of $\Lambda$, by taking maxima over all possible $\Cr{hermite}$ for $r \leq (nmd)^{2}$, one can assume $\Cr{hermite}$ to not depend on $r$ and only depend on $n,m,d$.

We also have the following important result that tells us that the Hadamard ratio is bounded from below. 
\begin{definition}
    \label{de:hadamard_ratio}

    Given a lattice $\Lambda \subseteq V$ in a Euclidean space $V$ with a $\mathbb{Z}$-basis $v_{1},\dots,v_{r}$. Then, the following quantity is called the Hadamard ratio of the basis $v_{1},\dots,v_{r}$.
    \begin{equation}
        \frac{\|v_{1}\| \|v_{2}\| \dots \|v_{r}\|}{H(\Lambda)} .
    \end{equation}
\end{definition}

If a basis $v_{1},\dots,v_{r}$ is an orthogonal basis, then it is clear that $\|v_{1}\| \dots \|v_{r}\| = H(\Lambda)$ and the Hadamard ratio is 1. For general bases, we only have the so-called Hadamard inequality: 

\begin{lemma}
    \label{le:hadamard_bound}
    Consider the same setup as \cref{de:hadamard_ratio}. Then,
    \begin{equation}
        \frac{\|v_{1}\| \|v_{2}\| \dots \|v_{r}\|}{H(\Lambda)} \geq 1.
    \end{equation}
    That is, the Hadamard ratio of a lattice is at least 1. 
\end{lemma}

\subsection{Hypothesis on test functions}

\label{se:test_funcs}

The functions that are of interest in this theory are compactly supported smooth functions and functions that are indicators of sets with nice boundaries. The following class of functions contains both of these cases.

\begin{hypothesis} We call a test function
    $f: \mathbb{R}^{d} \rightarrow \mathbb{R}$ 
    ``admissible'' if it is a 
    compactly supported bounded measurable function such that 
    the error function 
    \begin{equation}
        E_{f}(x, \varepsilon) = \sup_{\|x-y\| \leq \varepsilon} |f(x)-f(y)|
        \label{eq:defi_of_E}
    \end{equation}
    satisfies for some $\Cl[c]{admis} = \Cr{admis}(f)>0$, for every $1 \geq \varepsilon>0$ 
    and for any non-zero real subspace $V \subseteq \mathbb{R}^{d}$
    \begin{equation}
        \int_{V} E_{f}(x,\varepsilon) \diff x \leq \Cr{admis} \cdot \varepsilon,
    \end{equation}
    The integration is happening with the induced Lebesgue 
    measure from the inclusion $V \subseteq \mathbb{R}^{d}$ and $\Cr{admis}>0$ is required to be independent of $V$. 
    \label{hy:admissible}
\end{hypothesis}
\begin{remark}
    \label{re:help}
    \begin{enumerate}

        \

        \item One can replace the condition $0 < \varepsilon \leq 1$ by $0 < \varepsilon \leq \Cl[c]{someconst}$ for any constant $\Cr{someconst}>0$ by suitably updating the constant $\Cr{admis}(f)$. 
        \item 
            We observe that the hypothesis must assume that $\varepsilon$ is bounded, otherwise the only function satisfying the hypothesis would be $f=0$. 
            We thank Fabian Gundlach for pointing this out to us.
    \end{enumerate}
\end{remark}

A consequence of the above hypothesis is the following estimate for Riemann sums. 
\begin{lemma}
    \label{le:Riemann_estimate}

    Let $V$ be a Euclidean space and $\Lambda \subseteq V$ be a lattice such that $\dim \Lambda \otimes \mathbb{R} = n$.
    Let $f:V \rightarrow  \mathbb{R}$ be an admissible test function, in the sense of \cref{hy:admissible}.
    Let $\Cl[c]{arbit1} > 0$ be any arbitrary constant.
    Then, for some constant $\Cl[c]{admis-update}$ depending on $\Cr{arbit1}$ and the function $f$, for each
    $T \geq  \Cr{arbit1}^{-1} \rho(\Lambda)$ we have
    \begin{equation}
        \Big|\tfrac{1}{T^{n}}\sum_{v \in \Lambda} f(\tfrac{1}{T}v) - \tfrac{1}{H(\Lambda)} \int_{\Lambda \otimes \mathbb{R}}^{} f(x) \diff  x  \Big| \leq \frac{\Cr{admis-update} \cdot \rho(\Lambda)}{H(\Lambda)\cdot T},
    \end{equation}
    where the integral is with respect to the subspace measure on $\Lambda \otimes \mathbb{R} \subseteq V$. 
    Here 
    $\rho(\Lambda)$ denotes the covering radius of $\Lambda$.
\end{lemma}
\begin{proof}
    Let $F \subseteq V$ be a Voronoi domain of $\Lambda \subseteq \Lambda \otimes \mathbb{R}$.
    Then, by \cref{le:voronoi}, one gets $F \subseteq B_{0} (\rho)$ and that $\vol(F) = H(\Lambda)$.
    One then observes that
    {\small
        \begin{align}
            \Big|\tfrac{1}{T^{n}}\sum_{v \in \Lambda} f(\tfrac{1}{T}v) - \tfrac{1}{H(\Lambda)} \int_{\Lambda \otimes \mathbb{R}}^{} f(x) \diff  x  \Big|  
    &
    =  
    \Big|\tfrac{1}{T^{n}}\sum_{v \in \Lambda} f(\tfrac{1}{T}v) - \tfrac{1}{ T ^{n}H(\Lambda)  } 
    \int_{\Lambda \otimes \mathbb{R}}^{} f( \tfrac{1}{T} x) \diff  x  \Big|  
    \\ 
    &
    =  
    \tfrac{1}{H(\Lambda)T^{n}} \Big| \sum_{v \in \Lambda} f(\tfrac{1}{T}v) \int_{ F}^{} \diff x -    
    \sum_{v \in \Lambda}\int_{ F+v}^{}  f( \tfrac{1}{T} x) \diff  x  \Big|  
    \\ 
    &
    \leq  
    \tfrac{1}{H(\Lambda)T^{n}} \sum_{v \in \Lambda}\int_{ F+v}^{} \Big| f(\tfrac{1}{T}v) -    
    f( \tfrac{1}{T} x)   \Big| \diff x 
    \\ 
    & 
    \leq  
    \tfrac{1}{H(\Lambda)} \int_{x \in \Lambda \otimes \mathbb{R}}  E_f(x,\rho(\Lambda)/T) \diff x 
 \\ & \leq  
 \tfrac{1}{H(\Lambda)} \Cr{admis}(f) \frac{\rho(\Lambda)}{T}.
    \end{align}}
    Here $E_{f}(\cdot,\cdot)$ is as in \cref{eq:defi_of_E}. Note that for the last inequality obtained via Hypothesis \ref{hy:admissible} we set $\varepsilon = \rho(\Lambda)/T \leq \Cr{arbit1}$ by our assumed lower bound on $T$.
    The constant $\Cr{admis}(f)$ is updated to allow this choice of $\varepsilon >0$ (see~\cref{re:help})
    and $\Cr{admis-update}$ is simply chosen as $\Cr{admis}(f)$.
\end{proof}

\subsection{Summation by parts}

The following result on summation by parts will come in handy. It is a stronger form of~\cref{le:ballvol}.
\begin{lemma}
    \label{le:domain_dimension_bound}
    Let $\Lambda \subseteq V$ be a lattice in a Euclidean space $V$ as before.
    Let $h_1,h_2 \in \mathbb{Z}_{\geq 1}$. 
    Let $\mathcal{D} \subseteq V$ be a domain of infinite volume such that for $T \geq 1$ the domain $\mathcal{D}$ satisfies the following growth condition
    \begin{equation}
        \card\{v \in \Lambda  \cap \mathcal{D} \mid 
            \|v\| \leq T
        \}
        \leq 
        \Cl[c]{sbphyp}
        \cdot T^{h_1},
    \end{equation}
    for $\Cr{sbphyp} = \Cr{sbphyp}(\mathcal{D})$.

    Then for $1 \leq a \leq b$ one has:
    \begin{equation}
        \sum_{\substack{ l \in  \mathcal{D} \cap \Lambda  \\   a \leq \|l\| < b}} \frac{1}{\|l\|^{h_2}} 
        \leq \Cr{sbphyp}(\mathcal{D})\cdot
        \left( a^{h_1-h_2} + b^{h_1-h_2} + h_{2} \int_{a}^{b} x^{h_1 - h_2 -1} dx \right).
    \end{equation}
\end{lemma}
\begin{proof}
    As sketched in \cite[(13)]{K1994} we define 
    \begin{equation}
        \varphi(n) = \card \{ v \in \Lambda \cap \mathcal{D} : n\leq  \|v\| < n+1\}.
    \end{equation}
    Then one certainly has the inequality
    \begin{equation}
        \sum_{\substack{ l \in  \mathcal{D} \cap \Lambda  \\   a \leq \|l\| < b}} \frac{1}{\|l\|^{h_2}} 
 \leq \sum_{a\leq n<b}\varphi(n)\cdot n^{-h_2}
    \end{equation}
    Viewing $n^{-h_2}$ as a continuously differentiable function of $n$, Abel summation yields: 
    \begin{equation}
         \sum_{a\leq n< b}\varphi(n)\cdot n^{-h_2}\leq b^{-h_2}\cdot \Phi(b)- a^{-h_2}\cdot\Phi(a)+ h_2 \int_{a}^{b} x^{-h_{2}-1}\cdot\Phi(x) \diff x
    \end{equation}
    where $\Phi(x)=\sum_{n=0}^{\lceil x\rceil-1}\varphi(n)$ is the partial sum function. By assumption $\Phi(x)\leq\Cr{sbphyp}\cdot x^{h_1}$. The result follows.
\end{proof}

\subsection{Number fields, Euclidean structure and algebraic integer lattice}
\label{se:haar_on_KR}

Throughout this paper, we assume $K$ to be a number field of signature $(r_1,r_2)$, so that $r_1 + 2 r_2 = \deg K = d$. 
We fix, for once and for all, the following $l^{2}$-norm on $K \otimes \mathbb{R} \simeq \mathbb{R}^{r_1} \times \mathbb{C}^{r_2}$
\begin{equation}
    \|x\|^{2} = |\Delta_K|^{-\frac{1}{d}}  \Tr(x\overline{x}),
    \label{eq:norm}
\end{equation}
where $\Delta_K$ is the discriminant of the number field. The involution $\overline{(\ ) }$ in \eqref{eq:norm}
denotes complex conjugation
on all the complex places of $K$.
For any $r \in \mathbb{Z}_{\geq 1}$, the Euclidean space $K^{r} \otimes \mathbb{R} = (K\otimes \mathbb{R})^{r}$ comes equipped with the structure from $r$-fold copies of this underlying inner product. 
In particular, this also defines an $l^{2}$-norm on $\M_{n \times m}(\KR) \simeq \KR^{n \times m}$ for any $m,n > 0$.

It is known since the time of Minkowski that $\OK \subseteq \KR$ is a lattice, i.e. with respect to any Euclidean measure on $\KR$, $\vol(\KR/\OK) < \infty$.
Our quadratic form in \cref{eq:norm} is engineered to ensure that the lattice $\OK^{r} \subseteq \KR^{r}$ has unit covolume for any $r \geq 1$. 

We will fix a $\mathbb{Z}$-basis of $\OK$ for once and for all. Most of our implicit constants will depend on the choice of this basis. Here is a lemma demonstrating how this basis affects the underlying constants.
\begin{lemma}
    \label{le:basis_of_OK}
    Let $\Lambda \subseteq \KR^{m}$ be a free $\OK$-module of rank $k \leq m$, that is $ \Lambda = \OK v_{1} \oplus  \cdots \oplus \OK v_{k}$. 
    Then, there exists a $\mathbb{Z}$-basis $w_{11},\dots,w_{1d}, w_{21}, \dots, w_{(k-1)d}, w_{k1},\dots,w_{kd}$ such that for all $i=1,\dots,k$ of $\Lambda$
    and $j=1,\dots,d$ we have
    \begin{equation}
        \Cl[c]{okl}\|v_{i}\|\leq \|w_{ij}\| \leq \Cl[c]{okr} \|v_{i}\|,
        \label{eq:ij}
    \end{equation}
    where $\Cr{okl}$ and $\Cr{okr}$ depend on
    the $\mathbb{Z}$-basis of $\OK$ that we have fixed
    but neither on $v_{1},\dots,v_{k}$, nor on $\Lambda$ and not even on $k$ and $m$.
\end{lemma}
\begin{proof}
    Let $u_{1},\dots,u_{d}$ be our preselected $\mathbb{Z}$-basis
    so that $ \OK =  \mathbb{Z} u_{1} \oplus \cdots \oplus  \mathbb{Z}u_{d}$.
    Each $u_{i} \in \OK$ must clearly be non-zero.

    For $i,j$ as in \cref{eq:ij}, we choose
    $w_{ij} = u_{j} v_{i}$. 
    Then, clearly
    \begin{equation}
        \Big(\min_{\sigma: K \rightarrow  \mathbb{C}} \min_{j \in \{1,\dots,d\}}  |\sigma(u_{j})| \Big) \cdot \|v_{i}\|
        \leq
        \|w_{ij}\| 
        \leq \Big(\max_{\sigma: K \rightarrow  \mathbb{C}} \max_{j \in \{1,\dots,d\}}  |\sigma(u_{j})|  \Big)\cdot \|v_{i}\|.
    \end{equation}
\end{proof}

\section{Matrices, subspaces and lattices}

\label{se:schmidt}

\subsection{A useful trijection}
For an echelon matrix $D \in \M_{k \times m}$, let us define a lattice $\Lambda_{D}$ as follows.
\begin{definition}
    \label{de:lambda}
    Let $D \in \M_{k \times m}(K)$ be an echelon matrix. Then
    \begin{align}
        \Lambda_D &  = ( \M_{1 \times k}(K)\cdot D ) \cap \M_{1 \times m}(\OK).
        \label{eq:defi_of_lambda}
    \end{align}
    Hence $\Lambda_D$ is a lattice in $\M_{1 \times m}(\KR)$ 
    that contains all the vectors in $\M_{1 \times k}(K) \cdot D$ with integer entries. 
    It is an $\OK$-module. 
\end{definition}

We note that $\Lambda_{D}\subseteq \M_{1 \times m}(K)$ lives in a subspace of $K$-dimension $k < m$.
Observe that the following equalities hold.
\begin{align}
    \Lambda_D \otimes \mathbb{R} & = \Lambda_D \otimes_{\mathbb{Z}} \mathbb{R} = \Lambda_D \otimes _{\OK} K_{\mathbb{R}},\\
    \Lambda_D \otimes \mathbb{Q} & = \Lambda_D \otimes_{\mathbb{Z}} \mathbb{Q}  = \Lambda_D \otimes_{\OK } K .
\end{align}

To deal with our counting problems, we will need several equivalent descriptions of echelon matrices. The following proposition serves as a useful tool.
\begin{proposition}
    \label{pr:trijection}
    The following sets are in bijection with each other.
    \begin{enumerate}
        \item Rank $k$ row-reduced echelon matrices in $\M_{k \times m}(K)$.
        \item Points in $\Gr(k,K^{m})$.
        \item $\OK^{m}$-primitive $\OK$-modules of rank $k$ in $K^{m}$.
    \end{enumerate}
\end{proposition}

\begin{proof}
    The assignment $D \mapsto \Lambda_D \otimes \mathbb{Q}$ 
    assigns to an echelon matrix the rational subspace $ \M_{1 \times k}(K)D$ of $K$-dimension $k$ in $K^{m}$. One can recover $D$ from the subspace
    $\M_{1 \times k}(K) \cdot D$ by taking a $K$-basis and putting it in the appropriate echelon form.

    Note that the definition of $\Lambda_{D}$ in \cref{de:lambda} forces that $\Lambda_{D}$ is primitive. 
    The assignment $D \mapsto \Lambda_D$ is a bijection again since the echelon matrix $D$ can be recovered from $\Lambda_{D} \otimes \mathbb{Q}$ as described above.

    Here is a diagram showing this trijection with a third isomorphism laid out.

    \hspace{-0.45cm}
    \begin{tikzcd}
        && {\Gr(k,K^m)} \\
        \\
        \text{Echelon matrices in $\M_{k \times m}(K)$} && {\text{Rank-$k$ prim. } \OK\text{-modules in }\OK^m}
        \arrow["{D\mapsto \M_{1\times k}(K)\cdot D}"',leftrightarrow, from=1-3, to=3-1]
        \arrow["{S \mapsto S \cap \M_{1 \times m}(\OK)}",  leftrightarrow, from=1-3, to=3-3]
        \arrow["{D \mapsto \Lambda_D}",leftrightarrow, from=3-1, to=3-3]
    \end{tikzcd}

\end{proof}

Recall the definition of $H(~\cdot~)$ from \cref{de:height_definition}.
When there is no ambiguity, we shall at times write $H(\Lambda_D)$ as $H(V)$ or $H(D)$ for $V=\M_{1 \times k}(K) D$.  

\subsection{Schmidt's theorem}
The following is a result of W. Schmidt \cite{S1967}:
\begin{theorem}
    \label{th:schmidt}
    For $T \geq 1$, one has
    \begin{equation}
        \Cl[c]{schl}T^{m} \leq \card \{ V \in \Gr(k,K^{m}) \mid H(V) \leq T \} \leq \Cl[c]{schr} T^{m},
    \end{equation}
    for some constants $\Cr{schl},\Cr{schr}$ that depend on $K,k,m$.
\end{theorem}
In fact, more precise asymptotics were established by J. Thunder \cite{T1992}, but we will not require those.
The main point for us is that, using the trijection of \cref{pr:trijection}, there are finitely many echelon matrices $D\in \M_{k \times m}(K)$ satisfying $H(D) \leq T$.

A corollary of~\cref{th:schmidt} is the following lemma, which is also given in \cite[Corollary 17]{GSV23} but which we outline here for the sake of completeness.

\begin{lemma}
    \label{le:schmidt_makes_c1_finite}
    The constant $\Cr{main}$ defined in~\cref{eq:value_of_c} is finite for any admissible function $f: \M_{n \times m }(\KR) \rightarrow  \mathbb{R}$.
\end{lemma}
\begin{proof}
    This follows from the claim that for any echelon matrix $D \in \M_{k \times m}(K)$, one has 
    \begin{equation}
        \mathfrak{D}(D)^{-n}\int_{\M_{n \times k}(\KR)}^{}\Big| f(xD) \Big|\diff x  \leq \Cl[c]{height} H(D)^{-n},
        \label{eq:summable_matrices}
    \end{equation}
    for some constant $\Cr{height}$ depending on $f$. Once we have established this claim, we note that $n>m$ and \cref{th:schmidt} are sufficient to prove that the right-side of~\cref{eq:summable_matrices} is finitely summable over all echelon matrices in $\M_{k \times m}(K)$. 

    Indeed, there is a relation between $\mathfrak{D}(D)$ and $H(D)$. The product of the Jacobian of the map $x \rightarrow x D$ for $x \in M_{1 \times k}(K_\mathbb{R})$ times the factor $\mathfrak{D}(D)$ is exactly $H(D)$. 
   Indeed, $\mathfrak{D}(D)$ is the index by which the lattice $\Lambda_D$ differs from the image of $\mathcal{O}_{K}^{n}$ under $x\mapsto xD$, while the Jacobian accounts for the Euclidean volume distortion.
    This implies in particular that for any admissible function $f: \M_{n \times k}(\KR) \rightarrow  \mathbb{R}$ one gets (cf. Appendix A of \cite{GSV23})
    \begin{equation}
        \frac{1}{\mathfrak{D}(D)^{n}} \int_{\M_{n\times k}(K_\mathbb{R})} f( x D) \diff x = \int_{\M_{n \times k}(\KR) \cdot D} f(x) {\diff_{D}} x,
        \label{eq:d_D_defined}
    \end{equation}
    where $\mathrm{d}_D x$ is a Lebesgue measure on $\M_{n \times k}(\KR) \cdot D \subseteq \M_{n \times m}(\KR)$ 
    such that the lattice $\M_{n \times m}(\OK) \cap \M_{n \times k}(\KR) \cdot D  \subseteq \M_{n\times k}(K_\mathbb{R})\cdot D$
    has unit covolume. 
    See \cite[\S 3]{GSV23} for details. The measure $\diff _D x$ can then be expressed in terms of the induced Lebesgue measure $\diff_{l} x$ on $\M_{n\times k}(K_\mathbb{R})\cdot D \subseteq \M_{n \times m}(\KR)$ via the relation $\diff_{D} x  = H(D)^{-n} \diff_{l} x$. 
    We pick $\Cr{height}$ by setting
    \begin{equation}
        \Cr{height} = \sup_{\text{ Subspace }V \subseteq \M_{n \times m}(\KR)} \int_{V}^{} |f(x)|\diff_{l} x.
    \end{equation}
\end{proof}

\subsection{Successive minima in number fields}

\begin{definition}
    Consider $\KR^{m}$ as a Euclidean space equipped with the norm described by~\cref{eq:norm}. 
Endow $\KR^{m}$ with a total order by lexicographical ordering of coordinates via the fixed isomorphism $\KR^{m} \simeq \mathbb{R}^{d \times m}$.
    
Let $\Lambda \subseteq \KR^{m}$ be a lattice such that it is also an $\OK$-module of $\OK$-rank $k \leq m$. 
    Define the successive $K$-minima of $\Lambda$ as
    $l_i = l_i(\Lambda)$ for 
    $i=1,\dots,k$ given by 
    \begin{align}
        l_1(\Lambda) &  = \argmin_{v \in \Lambda \setminus \{0\} } \|v\| \\
        l_2(\Lambda)  & = \argmin_{v \in \Lambda \setminus K\cdot l_1  } \|v\| \\
        l_3(\Lambda)  & = \argmin_{v \in \Lambda \setminus K \cdot l_1 + K \cdot l_2  } \|v\| \\
                      &~\vdots
    \end{align}
    \label{de:defi_of_successive}
    Here $\argmin_{v \in A} \|v\|$ is the lexicographically least minimizer of the norm $\| \cdot\|$ among vectors in a non-empty discrete set $A \subseteq \KR^{m}$.
\end{definition}
\begin{remark}
The actual choice of the lexicographical ordering does not matter. The use of the total ordering on $\KR^{m}$ serves only to unambiguously define $\argmin$. We will make no use of this ordering in our arguments.
\end{remark}
Here is a lemma concerning the above definition.
\begin{lemma}
    \label{le:props_of_minima}
    Let $\Lambda \subseteq \KR^{m}$ be as in~\cref{de:defi_of_successive} of $\OK$-rank $k$
    and let $\{l_{i}(\Lambda)\}_{i=1}^{k}$ be the corresponding minima.
    Then, the following statements hold.
    \begin{enumerate}
        \item \label{lepart:1}
            Let $H(\Lambda)$ be the height of $\Lambda$ as defined in~\cref{de:height_definition}. 
            We have the relations
            \begin{equation}
                \|l_1\| \leq \Cr{hermite} H(\Lambda)^{\frac{1}{d k}}
                \text{ and }
                \|l_1\|^{d} \dots  \|l_{k}\|^{d} \leq {\Cl[c]{okhadamard} } H(\Lambda).
            \end{equation}
            Here the constants $\Cr{hermite}, \Cr{okhadamard}>0$ are independent of $\Lambda$.
        \item 
            \label{lepart:2}
            We get that 
            \begin{equation}
                \rho(\Lambda) \leq \Cl[c]{covering} \cdot \|l_{k}(\Lambda)\|
                ,
            \end{equation}
            where $\Cr{covering}>0$ is independent of $\Lambda$.
        \item \label{lepart:3}
            For $i<j$, denote by $\pi_{i}:\KR^{m} \rightarrow K_\mathbb{R} \cdot l_i$ 
             the orthogonal projection onto $K_\mathbb{R} \cdot l_i$. Then
            \begin{equation}
                \| \pi_{i}(l_j)\| \leq \Cl[c]{ijconst}
                \|l_i\|,
            \end{equation}
            where $\Cr{ijconst}$ does not depend on $i,j$ or $\Lambda$.
    \end{enumerate}

\end{lemma}
\begin{proof}

    { \em Part \ref{lepart:1}: }

    The first follows from~\cref{le:minkowski1} since $\|l_{1}(\Lambda)\|$ is the length of a shortest vector in $\Lambda$.
    The second statement follows from~\cite[Theorem 2]{FD10}.

    { \em Part \ref{lepart:2}: }

    We know that $\rho(\Lambda)$ is the smallest radius $r$ such that there exists some compact set $F \subseteq \Lambda \otimes \mathbb{R}$ such that 
    $F + \Lambda = \Lambda \otimes \mathbb{R}$ and $F \subseteq B_{0}(r)$. We will show that for some constant $\Cr{covering}>0$,
    $\Cr{covering} \|l_{k}\|$ works as such a radius $r$. For this, we will find a fundamental domain $F$ such that $F + \Lambda = \Lambda \otimes \mathbb{R}$ and $F \subseteq B_{0}(\Cr{covering} \|l_{k}\|)$.
Since we know that $\Lambda' = \OK l_1 + \dots + \OK l_k$ is a sublattice of finite index inside $\Lambda$, it is sufficient to find a domain $F$ such that $F + \Lambda' = \Lambda \otimes \mathbb{R}$ and $F \subseteq B_{0}(\Cr{covering} \|l_{k}\|)$.
    
For this, we consider a fundamental domain $F_{0} \subseteq \KR$ such that $F_{0} + \OK = \KR$. 
Observe that 
\begin{align}
     &  
    \sum_{i=1}^{k}(F_{0} + \OK) l_{i} = \sum_{i=1}^{k} \KR l_{i} = \Lambda \otimes \mathbb{R} \\ 
    \implies & 
    \sum_{i=1}^{k}F_{0} l_{i} + \Lambda'  = \Lambda \otimes \mathbb{R} 
\end{align}
So we can set $F = \sum_{i=1}^{k}F_{0} l_{i}$. One can then find a constant depending on $F_{0}$, which only depends on $\OK$, such that
\begin{equation}
v \in F_{0} x \implies  \|v\| \leq \Cl[c]{coverprim} \|x\| \text{ for all }x \in \KR \setminus\{0\}.
\end{equation}
Then it is clear that $\Cr{covering} = n \Cr{coverprim}$ works to imply that $F \subseteq B_{0}(\Cr{covering} \|l_{k}\|)$. Hence we are done.
    %

    { \em Part \ref{lepart:3}: }

    Observe that $l_j + \OK \cdot l_i \subseteq \Lambda$. We also know that for any $\alpha \in \OK$,
    the definition of $l_j$ implies $\|l_j\| \leq \|l _j + \alpha \cdot l_i\|$.
    It is clear that 
    \begin{equation}
        \pi_{i}(l_j + \alpha  \cdot l_i) = \pi_{i}(l_j) +   \pi_{i}(\alpha \cdot l_{i}) .
    \end{equation}
    Now $\alpha l_{i} \in \OK \cdot l_{i} \subseteq K_\mathbb{R} \cdot l_{i}$ so $\pi_{i}(\alpha \cdot l_i) = \alpha \cdot l_{i}$. Furthermore, we also know that for any $x \in M_{1 \times m}(K_\mathbb{R})$
    \begin{equation}
        \|x\|^{2} = \|\pi_{i}(x)\|^{2} + \|\pi_{i}^{\perp}(x)\|^{2},
    \end{equation}
    where $\pi_{i}^{\perp}$ is the projection to the orthogonal complement of $K_\mathbb{R}\cdot l_{i} \subseteq M_{1 \times m}(K_\mathbb{R})$. 
    We know that $\pi_{i}^{\perp}(l_j + \alpha \cdot l_i) = \pi_{i}^{\perp}(l_j)$.
    The net result is that 
    \begin{align}
& \|\pi_{i}(l_j) + \alpha \cdot l_i \|^{2} +  \|\pi_{i}^{\perp}(l_j)\|^{2}    = \|l_j + \alpha \cdot l_i\|^{2} \geq  \|l_j\|^{2} = \|\pi_{i}(l_j)\|^{2} + \|\pi_{i}^{\perp}(l_j)\|^{2}\\
\Rightarrow
& \|l_j + \alpha \cdot l_i\|^{2} - \|l_j\|^{2} = \|\pi_{i}(l_j) + \alpha \cdot l_i\|^{2} - \|\pi_{i}(l_{j})\|^{2} \geq 0.
    \end{align}
    This tells us that
    \begin{align}
  & \|\pi_{i}(l_j)\| \leq \min_{\alpha \in \OK} \|\pi_{i}(l_j) + \alpha \cdot l_i\| \\
        \Rightarrow & \|\pi_{i}(l_j)\| \leq \rho( \OK \cdot l_i).
    \end{align}
    It follows 
    from the proof of Part \ref{lepart:2} of the statement
    that the covering radius $\rho(\OK \cdot l_i) \leq \Cr{covering} \|l_i\|$. Hence, we are done.
\end{proof}

\subsection{Matrices with rows from a lattice}
Let us introduce a convenient notation for matrices containing rows taken from a particular lattice.
\begin{definition}
    \label{de:defi_of_M_t}
    For any subset $R \subseteq M_{1 \times m}(K_\mathbb{R})$, we denote by $M_{n}(R)$ the set of matrices in $M_{n \times m}(\KR)$ whose rows only contain elements of $R$.
\end{definition}

For $D \in M_{k \times m}(K)$ an echelon matrix,
observe that $\M_{n}(\Lambda_D \otimes \mathbb{Q}) = \M_{n \times k}(K) \cdot D$ 
and $\M_{n}(\Lambda_{D} \otimes \mathbb{R}) = \M_{n \times k}(K_\mathbb{R}) \cdot D$. 
However, $\M_n(\Lambda_D) \subseteq  \M_{n \times k}(\OK) \cdot D$ for an arbitrary echelon matrix $D \in \M_{k \times m}(K)$. In fact, 
\begin{equation}
    [\M_{n \times k}(\OK)  D : \M_{n}(\Lambda_{D})] = \mathfrak{D}(D)^{n},
    \label{eq:denominator_index}
\end{equation}
where $\mathfrak{D}(D)$ is defined in \cref{de:denonimator}.

Here is an important consequence of~\cref{eq:d_D_defined}.
\begin{lemma}
    \label{le:without_rank_cond}
    Suppose we have an admissible function $f: \M_{n \times m}(\KR) \rightarrow  \mathbb{R}$. 
    Let $\Cl[c]{arbit}>0$ be some arbitrary constant.
    Then, there exists a constant $\Cr{admis2} > 0$ depending on $f$,$n$ and $\Cr{arbit}$ 
    such that for any $T   \geq \Cr{arbit}^{-1}  \cdot \rho(\Lambda_{D})$ one has
    \begin{equation}
        \Big| \tfrac{1}{T^{knd}}\sum_{v \in \M_{n}(\Lambda_{D})}^{} f(\tfrac{1}{T}v) -  \mathfrak{D}(D)^{-n} \int_{\M_{n \times k}(\KR)}^{} f(xD)\diff x \Big| \leq \frac{ \Cl[c]{admis2} \cdot \rho(\Lambda_{D})}{T \cdot H({D})^{n}}
    \end{equation}
\end{lemma}
\begin{proof}
    The covolume of $M_n(\Lambda_D)$ is related to the covolume $H(D)$ of $\Lambda_D$ by the relation 
    $$H(M_n(\Lambda_D))=H(D)^n.$$
    To see this, one can be convinced by choosing a suitable basis for $\M_n(\Lambda_D)$ from a basis of $\Lambda_D$.

    The result then follows from the more general \cref{le:Riemann_estimate} 
    after we check that $\rho(\M_{n}(\Lambda_{D})) \leq \sqrt{n} \rho(\Lambda_{D})$. 
    We set $\Cr{admis2} = \sqrt{n }\cdot \Cr{admis-update}$ and while invoking ~\cref{le:Riemann_estimate}, we set $\Cr{arbit1} = \sqrt{n}\Cr{arbit}$.
\end{proof}

\begin{lemma}
    \label{le:new_bijection}
    One can rewrite the bijection in \cref{eq:bijection} as
    \begin{align}
 & \{A \in \M_{n \times m }(\OK) \mid \rank(A) = k \}   
 =   \bigsqcup_{ \substack{ D \in \M_{k \times m}(K) \\ D \text{ echelon }}}
 \{ A \in \M_{n}(\Lambda_{D}) \mid \rank A  = k\}
    \end{align}
    where $M_n(\Lambda_D)$ is as described in Definition \ref{de:defi_of_M_t}. 
\end{lemma}
\begin{proof}
    All one needs to check is that for $C \in \M_{n \times k}(K)$, the condition that $A = C \cdot D \in \M_{n \times m}(\OK)$ implies that the rows of $A$ must consist of elements of $\Lambda_{D}$ by definition of $\Lambda_{D}$ (see \cref{de:lambda}), and vice versa.
\end{proof}

\section{Integer matrices of fixed rank}

We will now begin collecting stepping stones towards establishing our main theorem.

\subsection{Matrices that interact with the support of the function}

First we introduce the following notation. 
For $1 \leq l \leq k$, define:
\begin{align}
  & \mathcal{F}_{l}(T)  = 
  \mathcal{F}_{l}^{(\Cr{sup})}(T) = 
  \\  &  \{D \in \M_{l \times m}(K), D \text{ is echelon, } \exists A \in \M_{n}(\Lambda_D),  \text{ $\rank A = l$ and } \|A\| \leq \Cl[c]{sup}T\},
  \label{eq:defi_of_calF}
\end{align}
where $\Cr{sup}$ is to be chosen later. 
For $l=0$, we define $\mathcal{F}_{0}(T) = \{0\}$.

The goal of defining $\mathcal{F}_{k}(T)$ is to identify matrices $D$ such that the sum  of $f$ given by $    \sum_{\substack{A \in \M_{n}(\Lambda_{D}) , \rank A = k }} f(\tfrac{1}{T} A )    $ is potentially non-zero.
The choice of $\Cr{sup}$ therefore has to be adjusted as per how large the support $\supp(f)$ of the function $f$ is.

\begin{lemma}
    \label{le:crude_early}
    One has for any choice of $\Cr{sup} > 0$ that 
    \begin{equation}
        D \in \mathcal{F}_{k}(T)  \implies  H(D) \leq \Cl[c]{crude2} T^{kd}.
    \end{equation}
    where $\Cr{crude2}>0$ depends on $\Cr{sup}$ and $\Cr{okl}$.
\end{lemma}
\begin{proof}
    Let $D \in \mathcal{F}_{k}(T)$. Consider some $A \in \M_{n}(\Lambda_{D})$
    with $\rank A  = k$. 
    Because $\rank A = k$, we know that the rows of $A$ contain a full-rank $K$-basis $\{v_{1},\dots,v_{k}\}$ of $\Lambda_{D} \otimes \mathbb{Q}$. 
    Let $\Lambda = \OK v_{1} \oplus \dots \oplus \OK v_{k}$ be the free $\OK$-module of $\OK$-rank $k$ generated by the basis vectors $v_{i} \in \M_{1 \times m}(K)$.
    Since any $\OK$-module $\Lambda$ generated by the rows of $A$ is a sublattice of $\Lambda_{D}$, we get $H(\Lambda) \geq H(\Lambda_{D})$. 

    Let us use the $\OK$-basis ${v_{i}}_{1 \leq i \leq k}$ of $\Lambda$ and obtain a $\mathbb{Z}$-basis $\{w_{ij}\}_{ {1 \leq i \leq k, 1 \leq j \leq d }}$
    of $\Lambda$ from \cref{le:basis_of_OK}. 
      After decreasing $\Cr{okl}$ if necessary, we may assume
  $\Cr{okl}\Cr{okr}\leq 1$.
    This tells us that 
    \begin{equation}
        \|A\|^{2} \geq \sum_{i=1}^{k} \|v_{i}\|^{2} \geq  \tfrac{\Cr{okl}^{2}}{d}\sum_{i=1}^{k}\sum_{j=1}^{d} \|w_{ij}\|^{2}.
    \end{equation}
    By \cref{le:hadamard_bound}, we know that the Hadamard ratio is bounded 
    from below for any $\mathbb{Z}$-basis $\{w_{ij}\}_{1 \leq i \leq k, 1 \leq j \leq d}$ of $\Lambda$.
    Then, the arithmetic-geometric means inequality gives us
    \begin{equation}
        \tfrac{1}{\Cr{okl}^{2} k }\|A\|^{2} \geq \tfrac{1}{kd}\sum_{i=1}^{k} \sum_{j=1}^{d} \|w_{ij}\|^{2} \geq  \Big( \prod_{i=1}^{k} \prod_{j=1}^{d} \|w_{ij}\|\Big)^{\frac{2}{kd}} \geq H(\Lambda)^{\frac{2}{kd}} \geq H(D)^{\frac{2}{kd}}.
    \end{equation}
    We now use the fact that $\|A\| \leq \Cr{sup} T$. Setting $\Cr{crude2} = ({\Cr{sup}/\sqrt{\Cr{okl}^{2} k}})^{kd} $ gives us the statement. 
\end{proof}
\begin{corollary}
    \label{co:counting_matrices}
    We have
    \begin{equation}
        \card \mathcal{F}_{k}(T) \leq  \Cl[c]{sizeoff} T^{mkd},
    \end{equation}
    where $\Cr{sizeoff} = \Cr{schr}\Cr{crude2}^{m}$.
\end{corollary}

We will also need the following lower bound on the height of matrices that are not in $\mathcal{F}_{k}(T)$.
\begin{lemma}
    \label{le:lower_bound_not_in_F}

    Let $D \in \M_{k \times m}(K)$ be an echelon matrix such that $D \notin \mathcal{F}_{k}(T)$.
    Then $H(D) \geq \Cl[c]{lowerF}T^{d}$ for some $\Cr{lowerF} > 0$ depending on $\Cr{sup}$.
\end{lemma}
\begin{proof}
    Let 
    $\{l_{1},\dots,l_{k}\} \subseteq \Lambda_{D}$ 
    be the successive minima of $\Lambda_{D}$ defined in~\cref{de:defi_of_successive}.
    By assumption on $D$, we must have that $\|l_{k}\| > k^{-\frac{1}{2}} \Cr{sup} T$ otherwise $D \in \mathcal{F}_{k}(T)$.
    Since $\Lambda_{D} \subseteq \M_{1 \times m}(\OK)$, there is a constant $\Cl[c]{minnorm}>0$, independent of $D$, such that $\|l_i\|\geq \Cr{minnorm}$ for every $i$.
    By~\cref{le:props_of_minima}, we also know that 
    \begin{equation}
        \|l_{1}\| \dots \|l_{k}\| \leq \Cr{okhadamard}^{1/d} H(D)^{1/d},
        \label{eq:also_useful}
    \end{equation}
    and therefore 
    \begin{equation}
        \Cr{minnorm}^{k-1} \cdot \Cr{sup} \cdot T \leq \Cr{okhadamard}^{1/d}H(D)^{1/d},
    \end{equation}
    and we are done.
\end{proof}

Using~\cref{le:lower_bound_not_in_F}, one gets the following convergence estimate for~\cref{le:schmidt_makes_c1_finite}.
\begin{corollary}
    \label{co:tail}

    We have 
    \begin{align}
  &  \sum_{\substack{D \in \M_{k \times m}(K),~D \text{ echelon } \\ D \notin \mathcal{F}_{k}(T)}}^{}\mathfrak{D}(D)^{-n} \int_{x \in \M_{n \times k}(\KR)}^{} | f(xD) |\diff x   \\
  & \leq
   \sum_{\substack{D \in \M_{k \times m}(K),~D \text{ echelon } \\ H(D) \geq \Cr{lowerF} T^{d}}}^{}\mathfrak{D}(D)^{-n} \int_{x \in \M_{n \times k}(\KR)}^{} |f(xD) | \diff x \leq  \Cl[c]{tailconst}  \tfrac{1}{T^{d}}
    \end{align}
    Here $\Cr{tailconst}>0$ is a constant that does not depend on $T>0$.
\end{corollary}
\begin{proof}
  As in the proof of~\cref{le:schmidt_makes_c1_finite}, one has
  \[
  \mathfrak{D}(D)^{-n}
  \int_{\M_{n \times k}(\KR)} |f(xD)|\diff x
  \leq \Cr{height} H(D)^{-n}.
  \]
  Hence, by~\cref{le:lower_bound_not_in_F}, the left-hand side is bounded by
  \[
  \Cr{height}
  \sum_{\substack{D \in \M_{k \times m}(K),~D\text{ echelon}\\
  H(D)\geq \Cr{lowerF}T^{d}}}
  H(D)^{-n}.
  \]
  Summation by parts and~\cref{th:schmidt} show that this is at most
  \[
  \Cr{tailconst}T^{-d(n-m)}
  \leq \Cr{tailconst}T^{-d},
  \]
  since \(n-m\geq1\).
\end{proof}

Using the notation $\mathcal{F}_{k}(T)$, we can rewrite the sum in~\cref{eq:main_eq}: 
\begin{lemma}
    \label{le:crude_bound}
    The left-hand side of \cref{eq:main_eq} satisfies
    \begin{equation}
        \label{eq:left_side}
        \sum_{\substack{A \in \M_{n \times m}(\OK) \\ \rank A = k}} f( \tfrac{1}{T} A) =
        \sum_{D \in \mathcal{F}^{(\Cr{sup})}_{k}(T)}^{}  \sum_{\substack{A \in \M_{n}(\Lambda_{D}) \\ \rank A = k}} f( \tfrac{1}{T} A),
    \end{equation}
    for some $\Cr{sup} > 0$ in~\cref{eq:defi_of_calF}.
\end{lemma}
\begin{proof}
    Indeed, one can set
    \begin{equation}
        \Cr{sup} = 1 + \sup\{ \|A\| : A \in \supp(f)\}.
        \label{eq:choice_of_sup}
    \end{equation}
    The bijection in~\cref{le:new_bijection} then allows one to conclude~\cref{eq:left_side}. The echelon matrices $D \in \M_{k \times m}(K)$ that are not in $\mathcal{F}_{k}(T)$ do not contribute to the sum due to the choice of $\Cr{sup}$.
\end{proof}

From now on, whenever we mention $\mathcal{F}_{k}(T)$, 
we assume that $\Cr{sup}$ has been chosen so that~\cref{le:crude_bound} holds. 
Note that this choice, given in~\cref{eq:choice_of_sup}, does not depend on $k$.
We will eventually use \cref{le:crude_bound} to prove \cref{th:main} in \cref{ss:proof_of_main}.

\subsection{Possible successive minima}

One has the following correspondence between matrices in $\mathcal{F}_{k}(T)$ which will be used in the proof of~\cref{th:main}.
\begin{lemma}
    \label{le:injective_minima}

    Denote for a constant $\Cl[c]{minima}>0$ the set
    \begin{align}
        \mathcal{B}_{k}^{(\Cr{minima})} (T)= \{(l_{1},\dots,l_{k}) \in \M_{1 \times m}(\OK)^{k} \mid  \substack{\|l_{1}\| \leq \dots \leq \|l_{k}\| \leq \Cr{minima} T \\  \|\pi_{i}(l_{j})\| \leq \Cr{minima} \|l_{i}\|  \text{ for each }j>i }\},
    \end{align}
    where the map $\pi_{i}:\M_{1 \times m}(\KR) \rightarrow  \M_{1 \times m}(\KR)$ is orthogonal projection onto $\KR \cdot l_{i}$.

    For each $D \in \mathcal{F}_{k}(T)$, consider the correspondence $D \mapsto  \{l_{i}(\Lambda_{D})\}_{i=1}^{k}$, where $l_{i}(\Lambda_{D})$ are the successive $K$-minima defined in~\cref{de:defi_of_successive}. Then, 
    \begin{enumerate}
        \item 
            For some $\Cr{minima} > 0$ that does not depend on $T$, the image of the map lies in $\mathcal{B}_{k}^{(\Cr{minima})}(T)$ for every $T$.
        \item This mapping is injective.
    \end{enumerate}
\end{lemma}
\begin{proof}
    The only thing to check is that the successive minima satisfy the properties demanded by $\mathcal{B}_{k}^{(\Cr{minima})}(T)$.
    We will use \cref{le:props_of_minima}.
    Clearly, the property $\|\pi_{i}(l_{j}) \|\leq \Cr{minima} \|l_{i}\|$ holds due to the third part of ~\cref{le:props_of_minima} for $\Cr{minima}$ chosen appropriately. 

    To get $\|l_{k}\| \leq \Cr{minima} T$, we note that $D \in \mathcal{F}_{k}(T)$. Therefore we know that the set $\{v \in  \Lambda_{D} \mid \|v\| \leq \Cr{sup}T \}$ contains at least $k$ vectors that are $K$-linearly independent. Indeed, this follows from how $\Cr{sup}$ is defined in~\cref{eq:defi_of_calF}. This implies in particular that $\|l_{k}\| \leq \Cr{sup} T$, since $K \cdot l_{1} + \dots + K\cdot  l_{k-1}$ cannot cover all the points in the set $\{ v \in \Lambda_{D} \mid \|v\| \leq \Cr{sup} T \}$. 
    Hence up to replacing $\Cr{minima}$, one can ensure that $\Cr{minima} \geq \max\{\Cr{ijconst} ,\Cr{sup}\}$ and we get $\|l_{1}\| \leq \dots \leq \|l_{k}\| \leq  \Cr{minima} T$.

    For the second part, one can recover $\Lambda_{D}$ from $\Lambda_{D} \otimes \mathbb{Q}$ due to the trijection in~\cref{pr:trijection}. This concludes the proof.
\end{proof}

When we will invoke~\cref{le:injective_minima} in the proof of~\cref{th:main}, we will assume that $\Cr{minima}>0$ has been chosen large enough for the conclusion of ~\cref{le:injective_minima} to be true. Then, we will refer to $\mathcal{B}_{k}^{(\Cr{minima})}(T)$ as $\mathcal{B}_{k}(T)$.

\begin{corollary}
    \label{co:covering_bound}
    For $D \in \mathcal{F}_{k}(T)$, one has that 
    \begin{equation}
        \rho(\Lambda_{D}) \leq (\Cr{covering} \cdot \Cr{minima} )\cdot T,
    \end{equation}
    where $\Cr{minima}>0$ is large enough so that the conclusions of \cref{le:injective_minima} are true.
\end{corollary}
\begin{proof}
    Follows immediately from $\rho(\Lambda_{D}) \leq \Cr{covering} \| l_{k}(D)\|$.
\end{proof}

\subsection{Low rank terms}



We begin by the following lemma.
\begin{lemma}
    \label{le:low_rank_induction}

    Let $1 \leq l < k$ and $\mathcal{F}_{l}(T), \mathcal{F}_{k}(T)$ be as in~\cref{eq:defi_of_calF}. 
    Let $D' \in \mathcal{F}_{l}{(T)}$.
    For some $\Cl[c]{lrind} > 0$ which 
    depends on 
    $\rho(\OK)$ but not on $l,T$ or $D'$,
    one then has for all $T \geq 1$:
    \begin{equation}
        \card \{D \in \mathcal{F}_{k}(T) \mid \Lambda_{D'} \subseteq \Lambda_{D} \} \leq \Cr{lrind}T^{d(k-l)(m-l)} H(D')^{k-l}
    \end{equation}
\end{lemma}
\begin{proof}
    Let $D_1,D_2\in \mathcal{F}_{k}(T)$ such that $\Lambda_{D'} \subseteq \Lambda_{D_i}$ for $i=1,2$. 
    We know that each $\Lambda_{D_{i}}$ and $\Lambda_{D'}$ is a subset of $\M_{1 \times m}(\OK)$
    as per~\cref{eq:defi_of_lambda}.
    Then the set of vectors in $S=\{v \in \M_{1 \times m}(\OK) \mid  \|v\| \leq \Cr{sup} T\}$
    contain a $K$-basis of $\Lambda_{D_i} \otimes \mathbb{Q}$ for $i=1,2$. 
    Moreover, since $D' \in \mathcal{F}_{l}(T)$, we can also conclude that $S$ contains a $K$-basis of $\Lambda_{D'} \otimes \mathbb{Q}$.

    In particular, there exist primitive vectors $(l_j^{(i)})_{j=1}^{k-l}$ in $S \subseteq \M_{1 \times m}(\OK)$, 
    with $\lVert l_j^{(i)}\rVert\leq  \Cr{sup} T $ for all $i=1,2$ and $j=1,\dots,k-l$, such that
    $$\Lambda_{D_i}\otimes\Q=(\Lambda_{D'}\otimes\Q)~ \oplus ~\bigoplus_{j=1}^{k-l}\ (l_j^{(i)}\cdot K), $$for $i=1,2$. 
    Observe that $\Lambda_{D_1}\otimes\Q=\Lambda_{D_2}\otimes\Q$ if and only if the two $K$-spaces spanned by $(l_j^{(1)})_{j=1}^{k-l}$ and $(l_j^{(2)})_{j=1}^{k-l}$,  respectively, are equal modulo $\Lambda_{D'}\otimes\Q$.

    We would therefore like to bound the number of choices for each $l_i$ up to $\Lambda_{D'}$-equivalence. To that end, we bound the number of lattice points in $S$ after projection of $\M_{1 \times m}(\OK)$ 
    onto $(\Lambda_{D'} \otimes \mathbb{R})^\perp$. 
    This is a $\mathbb{Z}$-lattice of rank $d(m-l)$ and of height $H(D')^{-1}$.
    So from \cref{le:ballvol}, the number of choices for each $l_i$
    inside this projection is bounded by $$ \Cr{ballvol}\left(T+ \rho(\M_{1 \times m}(\OK))  \right)^{d(m-l)} \times \frac{1}{H(D')^{-1}}.$$

    Since we are choosing $k-l$ vectors we arrive at the upper bound in the statement.
\end{proof}

We will now use the following to bound the low-rank terms that do not appear on the left-hand side of~\cref{eq:sums_distributes}, but will be added artificially in the proof of~\cref{th:main} in~\cref{ss:proof_of_main}.
\begin{lemma}
    \label{le:low_rank_terms}

    Let $f: \M_{n \times m}(\KR) \rightarrow  \mathbb{R}$ be an admissible function. Then, for $T \geq 1$ we have 
    \begin{equation}
        \tfrac{1}{T^{knd}}    \sum_{ D \in \mathcal{F}_{k}(T)} \sum_{\substack{A \in \M_{n}(\Lambda_{D}) \\ \rank A < k}} |f(\tfrac{1}{T} A)|   \leq 
        \Cl[c]{lowrank} \frac{ 1+ \log T}{T^{d(n-m+k-1)}} 
        \label{eq:loewrank}
    \end{equation}
    for some constant $\Cr{lowrank}>0$. 
\end{lemma}
\begin{proof}
    Without loss of generality, to prove the upper bound we can replace $f$ by $|f|$ and assume that $f \geq 0$.
    %
    One has
    \begin{align}
        \tfrac{1}{T^{knd}} \sum_{D \in \mathcal{F}_{k}{(T)}} \sum_{\substack{A \in \M_{n}(\Lambda_{D}) \\ \rank A < k }}^{}  
        f(\tfrac{1}{T} A )
        = 
        \tfrac{1}{T^{knd}} \sum_{l=0}^{k-1} \sum_{D' \in \mathcal{F}_{l}{(T)}}^{} \sum_{\substack{A \in \M_{n}(\Lambda_{D'}) \\ \rank A = l }}^{}  
        f(\tfrac{1}{T}A) n_{k}(D'),
    \end{align}
    where 
    \begin{equation}
        n_{k}(D') = \card\{D \in \mathcal{F}_{k}{(T)}  \mid \Lambda_{D'} \subseteq \Lambda_{D} \}.
    \end{equation}
    Indeed, the equality can be derived by observing that for any $D \in \mathcal{F}_{k}(T)$ and $A \in \M_{n}(\Lambda_{D})$ contributing to the sum on the left, the rows of $A$ must all live in $\Lambda_{D'}$ for some matrix $D' \in \M_{l \times m}(\Lambda_{D})$ of strictly smaller rank $l$. Since $A$ is in the support of $f(\tfrac{1}{T} \cdot x )$, this means that $D' \in \mathcal{F}_{l}(T)$. Finally, the term $n_{k}(D')$ is exactly what is needed to correct for overcounting.

    Let us first handle the $l=0$ term. One has that 
\begin{align}
        \tfrac{1}{T^{knd}}  \sum_{D' \in \mathcal{F}_{l}{(T)}}^{} \sum_{\substack{A \in \M_{n}(\Lambda_{D'}) \\ \rank A = 0 }}^{}   
        f(\tfrac{1}{T}A) n_{k}(D') 
         & = 
        \tfrac{1}{T^{knd}}  f(0 ) \card \mathcal{F}_{k}(T)  \\ 
     &   \leq \Cr{sizeoff} f(0) T^{-(n-m)kd}.
\end{align}
From the rest of the proof, we will assume $l \geq 1$.

    By~\cref{le:low_rank_induction}, 
    we write that
    \begin{align}
      & \tfrac{1}{T^{knd}} \sum_{D \in \mathcal{F}_{k}{(T)}} \sum_{\substack{A \in \M_{n}(\Lambda_{D}) \\ 0 <  \rank A < k }}^{}  
      f(\tfrac{1}{T} A )
      \\
        \leq &
        \tfrac{\Cr{lrind}}{T^{knd}} \sum_{l=0}^{k-1} 
        \sum_{D' \in \mathcal{F}_{l}{(T)}}^{} 
        T^{d(k-l)(m-l)}H(D')^{k-l}
        \sum_{\substack{A \in \M_{n}(\Lambda_{D'}) \\ \rank A = l }}^{}  
        f(\tfrac{1}{T}A) 
        .
    \end{align}
    For $T \geq 1$, using~\cref{le:ballvol} and~\cref{eq:choice_of_sup},
    the innermost sum above can be bounded as
    \begin{align}
        \sum_{\substack{A \in \M_{n}(\Lambda_{D'}) \\ \rank A = l }}^{}  
        f(\tfrac{1}{T}A)  \leq
        \sum_{\substack{A \in \M_{n}(\Lambda_{D'}) \\ \|A\| \leq \Cr{sup} T } }^{} \|f\|_{\infty}
      & \leq \Cr{ballvol} \|f\|_{\infty} (\Cr{sup} T + \rho(\Lambda_{D'      }))^{lnd} H(D')^{-n} \\
      & \leq \Cl[c]{lowrank2} H(D')^{-n} T^{lnd}.
    \end{align}
    Note here that $\rho(\Lambda_{D}') \leq (\Cr{covering} \Cr{minima}) T$ from~\cref{co:covering_bound}.
    Here $\Cr{lowrank2}>0$ is taken to be at least $ \|f\|_{\infty} \cdot \Cr{ballvol} \cdot (\Cr{sup} + \Cr{covering} \Cr{minima})^{(k-1)nd}$.

    Now take $\Cl[c]{lowrank3} = \Cr{lowrank2} \cdot \Cr{lrind}$. We obtain
    \begin{align}
        \tfrac{1}{T^{knd}} \sum_{D \in \mathcal{F}_{k}{(T)}} \sum_{\substack{A \in \M_{n}(\Lambda_{D}) \\ 0  < \rank A < k }}^{}  
      & f(\tfrac{1}{T} A ) 
      \leq  \Cr{lowrank3} \sum_{l=1}^{k-1} T^{d((k-l)(m-l) + n(l - k))} 
      \sum_{D' \in \mathcal{F}_{l}(T)}^{}
      H(D')^{k-l-n}.
      \label{eq:inequality_final_desire}
    \end{align}
    We shall now use summation by parts to settle the main claim of the lemma. 
    One has from~\cref{th:schmidt},~\cref{le:crude_early} and summation by parts that
    \begin{align}
        \sum_{D' \in \mathcal{F}_{l}(T)}^{}
        H(D')^{k-l-n} 
      & \leq 
      \sum_{\substack{D' \in \M_{l \times m}(K) \\ D' \text{  is echelon}, H(D') \leq \Cr{crude2}T^{ld} }} H(D')^{k-l-n} 
      \\
      & \leq 
      {(\Cr{crude2}T^{ld})^{k-l-n}}
      \eta( \Cr{crude2}T^{ld})
      \\
      &\quad + (n+l-k)\int_{1}^{\Cr{crude2}T^{ld}} \eta(x) x^{k-l-n-1} \diff x   \\
      \label{eq:ineqluaty_abel}
    \end{align}
    where we have
    \begin{equation}
        \eta(x) = 
        \sum_{\substack{D' \in \M_{l \times m}(K) \\ D' \text{  is echelon}, H(D') \leq x }} 1 \ \ \ \ \ \   \leq \Cr{schr} \cdot x^{{m}}.
        \label{eq:ineqluaty_schmidt}
    \end{equation}
    The inequality in~\cref{eq:ineqluaty_schmidt} is due to~\cref{th:schmidt}.
    Putting~\cref{eq:ineqluaty_schmidt} in~\cref{eq:ineqluaty_abel} gives us that
    \begin{equation}
        \sum_{D' \in \mathcal{F}_{l}(T)}^{}
        H(D')^{k-l-n} 
        \leq \Cl[c]{crudenew}\cdot 
        \begin{cases}
          T^{ld(m+k-l-n)} & \text{ if } m+k-l-n > 0 \\
           1 + \log T&  \text{ if } m+k-l-n=0 \\ 
           1 &\text{ if }m+k-l-n < 0 
        \end{cases}
        \label{eq:inequality_semifinal}.
    \end{equation}
    Set $\alpha_{l} = d \Big( (k-l)(m-l)  + n(l-k)\Big)$ and set $B_{l}(T)$ to be the right hand side of~\cref{eq:inequality_semifinal}.
    Putting~\cref{eq:inequality_semifinal} in~\cref{eq:inequality_final_desire} gives us
    \begin{align}
        \tfrac{1}{T^{knd}} \sum_{D \in \mathcal{F}_{k}{(T)}} \sum_{\substack{A \in \M_{n}(\Lambda_{D}) \\ 0 < \rank A < k }}^{}  
      & f(\tfrac{1}{T} A )  \\
      & \leq  \Cr{lowrank3} \cdot \sum_{l=1}^{k-1} 
      T^{ \alpha_{l}} B_{l}(T)
    \end{align}

When $m+k-l-n > 0$, we know that $T^{\alpha_{l}} T^{ld(m+k-l-n)} = T^{d N }$ 
    where
    \begin{align}
        N  & =  (k-l)(m-l)  + n(l-k) + l(m+k-l-n)  \\
           & =    (l-k)(n-m+l) -l(n-m+l) +lk   \\
           & =    -k(n-m+l)  +lk   = -(n-m)k  \leq -1.
    \end{align}
    Hence we have negative powers of $T$. 

For the case when $m+k-l-n < 0$, one gets
    \[
        T^{\alpha_{l}} (1 + \log T) \leq T^{-d(n-m+k-1)} (1 + \log T),
    \]
    since
    \[
        d^{-1}\alpha_{l} = (m-l-n)(k-l) \leq -(n-m+k-1),
    \]
    with equality at $l=k-1$.
    Finally, when $m+k-l-n = 0$ one has $\alpha_{l}\leq 0$ so we can bound the sum by a constant. In any case, $\sum_{l=1}^{k-1} T^{\alpha_{l}} B_{l}(T)$ is at most a constant times $(1 + \log T)T^{-d(n-m+k-1)}$. Hence we are done.
\end{proof}

\subsection{Putting it all together}
\label{ss:proof_of_main}

\begin{proof} {\bf (of~\cref{th:main})}
    In order to show~\cref{eq:main_eq} with $\Cr{main}$ given in~\cref{eq:value_of_c}, it is enough to show that in the non-exceptional case when $(n-m)dk > 1$
    \begin{equation}
        \sum_{D \in \mathcal{F}_{k}(T)}^{}  \Big|
        \tfrac{1}{T^{knd}}\sum_{\substack{A \in \M_{n}(\Lambda_{D}) \\ \rank A = k}} f( \tfrac{1}{T} A)  - 
        \mathfrak{D}(D)^{-n}\int_{\M_{n \times k}(\KR)} f( x D )\, dx \Big| \leq \tfrac{1}{T}\Cl[c]{error2} ,
        \label{eq:required_eq}
    \end{equation}
    for some $\Cr{error2}>0$ as specified, otherwise with an additional $\log T$ factor in the exceptional case.

    Indeed, this is because the terms
    \begin{align}
    & \sum_{D \notin \mathcal{F}_{k}(T)}^{}  \Big|
    \tfrac{1}{T^{knd}}\sum_{\substack{A \in \M_{n}(\Lambda_{D}) \\ \rank A = k}} f( \tfrac{1}{T} A)  - 
    \mathfrak{D}(D)^{-n}\int_{\M_{n \times k}(\KR)} f( x D )\, dx \Big|    \\
        = & \sum_{D \notin \mathcal{F}_{k}(T)}^{}  \Big| 
        \mathfrak{D}(D)^{-n}\int_{\M_{n \times k}(\KR)} f( x D )\, dx \Big|  \\ 
        \leq & 
        \Cr{tailconst} \tfrac{1}{T^{d}},
        \label{eq:nonrequired_eq}
    \end{align}
    where we used~\cref{le:crude_bound} in the first step
    and~\cref{co:tail} for the final inequality.

    If one could drop the $\rank A= k$ condition from the sum over $\M_{n}(\Lambda_{D})$, then one could invoke~\cref{le:without_rank_cond}. 
    From~\cref{co:covering_bound}, the condition $T  \geq \Cr{arbit}^{-1}  \cdot \rho(\Lambda_{D})$ is satisfied 
    if we choose $\Cr{arbit} = \Cr{covering}\cdot \Cr{minima}$.
    Then, we can make some progress on proving~\cref{eq:required_eq} by writing
    \begin{align}
  & 
  \sum_{D \in \mathcal{F}_{k}(T)}^{}  \Big|
  \tfrac{1}{T^{knd}}\sum_{\substack{A \in \M_{n}(\Lambda_{D}) \\ \rank A = k}} f( \tfrac{1}{T} A)  - 
  \mathfrak{D}(D)^{-n}\int_{\M_{n \times k}(\KR)} f( x D )\, dx \Big|  \\
        \leq & 
        \sum_{D \in \mathcal{F}_{k}(T)}^{}  \Big|
        \tfrac{1}{T^{knd}}\sum_{\substack{A \in \M_{n}(\Lambda_{D}) }} f( \tfrac{1}{T} A)  - 
        \mathfrak{D}(D)^{-n}\int_{\M_{n \times k}(\KR)} f( x D )\, dx \Big|  \\
             & + \tfrac{1}{T^{knd}}\sum_{D \in \mathcal{F}_{k}(T)}  \sum_{\substack{A \in \M_{n}(\Lambda_{D}) \\ \rank A < k}}^{} | f(\tfrac{1}{T} A ) |.
    \end{align}
    From~\cref{le:low_rank_terms}, the last term is smaller than $\Cr{lowrank}T^{-d(n-m+k-1)}(1+ \log T)$, so it absorbs in the constant $\Cr{error2}$ without any problems.

    As explained before, we know that from~\cref{co:covering_bound} that
    $D \in \mathcal{F}_{k}(T) \implies  \rho(\Lambda_{D}) \leq (\Cr{covering} \cdot \Cr{minima} )\cdot T$. Hence from \cref{le:without_rank_cond} 
    after setting $\Cr{arbit} = \Cr{covering}\cdot \Cr{minima}$
    one has that 
    \begin{align}
        \sum_{D \in \mathcal{F}_{k}(T)}
        \Big|
        \tfrac{1}{T^{knd}}\sum_{\substack{A \in \M_{n}(\Lambda_{D}) }} f( \tfrac{1}{T} A)  -  
    & \mathfrak{D}(D)^{-n}\int_{\M_{n \times k}(\KR)} f( x D )\, dx \Big|   \\ 
    & \leq  \frac{\Cr{admis2}}{T}\sum_{D \in \mathcal{F}_{k}(T)}^{} 
    \frac{ \rho(\Lambda_{D})}{H(D)^{n}}.
    \label{eq:ineq212}
    \end{align}



    {\bf When $n=m+1$,$k=1$ and $d=1$:}
    
    We treat this case separately. We are looking at the sum
    \begin{equation}
        \sum_{\substack{A \in \M_{n \times m}(\mathbb{Z}) \\ \rank A = 1}} f(\tfrac{1}{T} A) 
        =  
        \tfrac{1}{2}\sum_{v \in \mathbb{Z}^{m}_{\prim}}^{} \sum_{w \in \mathbb{Z}^{n} \setminus \{0\} } f(\tfrac{1}{T} wv^{T}).
    \end{equation}
    This has been discussed in~\cref{ss:log_term}.
    
    {\bf When $(n-m)d > 1$:}
    The goal is to now show that 
    \begin{equation}
        \sum_{D \in \mathcal{F}_{k}(T)} \frac{\rho(\Lambda_{D})}{H(D)^{n}} \leq \Cl[c]{finitesum},
    \end{equation}
    for some $\Cr{finitesum} > 0$ that does not depend on $T$.

    Recall $\mathcal{B}_{k}(T)$ from \cref{le:injective_minima}. Using~\cref{le:props_of_minima} to get that
    $\rho(\Lambda_{D}) \leq \Cr{covering} \|l_{k}\|$ and $\|l_{1}\|^{d}\dots \|l_{k}\|^{d} \leq \Cr{okhadamard} H(D)$, 
    one gets that for $\Cl[c]{innervar}= \Cr{covering} \cdot \Cr{okhadamard}^{n}$ one has
    \begin{equation}
    \sum_{D \in \mathcal{F}_{k}(T)} \frac{\rho(\Lambda_{D})}{H(D)^{n}} 
    \leq 
    \Cr{innervar}
    \sum_{(l_{1},\dots,l_{k}) \in \mathcal{B}_{k}(T)}^{} \frac{1}{\|l_{1}\|^{nd} \dots \|l_{k}\|^{nd-1}}.
    \label{eq:just_as_before}
    \end{equation}
    When $(n-m)d > 1$, one can relax the extra conditions on $\mathcal{B}_{k}(T)$ and bound the sum by 
    \begin{equation}
    \sum_{\substack{(l_{1},\dots,l_{k}) \in  \M_{1 \times m}(\OK) \\  \|l_{1}\|  \leq \dots \leq \|l_{k}\| \leq  \Cr{minima} T  }   }^{} \frac{1}{\|l_{1}\|^{nd} \dots \|l_{k}\|^{nd-1}}.
    \end{equation}
    When $k \geq 2$, the innermost sum is 
    \begin{equation}
    \sum_{\substack{ l_k\in  \M_{1 \times m}(\OK) \\  \|l_{k-1}\|  \leq \|l_{k}\| \leq  \Cr{minima} T  }   }^{} \frac{1}{\|l_{k}\|^{nd-1}} \leq 
    \sum_{\substack{ l_k\in  \M_{1 \times m}(\OK) \\  \|l_{k-1}\|  \leq \|l_{k}\| < \infty}   }^{} \frac{1}{\|l_{k}\|^{nd-1}} \leq \Cl[c]{easycase} \|l_{k-1}\|^{-(n-m)d+1},
    \end{equation}
    Here in the last sum we used~\cref{le:domain_dimension_bound}. 
    We can proceed with the rest of the nested sums inductively. For the base case, when $k=1$ it is clear that the following sum converges
    \begin{equation}
    \sum_{\substack{ l_1\in  \M_{1 \times m}(\OK) \\  \|l_{1}\| \leq  \Cr{minima} T  }   }^{} \frac{1}{\|l_{1}\|^{nd-1}} < \infty.
    \end{equation}

    {\bf When $n=m+1$, $d=1$ and $k \geq 2$:}
This is the only case remaining and it requires a separate argument. 
We can now write $\mathcal{O}_{K} = \mathbb{Z}$ since $d=1$.
The goal is to bound again the sum in~\cref{eq:just_as_before} but with $d=1$.
Fortunately, this reduces exactly to the case that was already covered by~\cite[\S 4]{K1994}.

However, we will still write an argument here for completeness. Consider the critical innermost sum of
\begin{equation}
\sum_{(l_{1},\dots,l_{k}) \in \mathcal{B}_{k}(T)}^{} \frac{1}{\|l_{1}\|^{n} \dots \|l_{k-1}\|^{n}\cdot\|l_{k}\|^{n-1}}.
\end{equation}
Dropping some conditions on $l_{k}$ it is enough to show that the sum satisfies
\begin{equation}
\sum_{ \substack{l_{k} \in \M_{1 \times m}(\mathbb{Z}) \\ \|l_{k-1}\| \leq \|l_{k}\| \leq \Cr{minima} T  \\ \|\pi_{k-1}(l_{k})\| \leq \Cr{minima} \|l_{k-1}\| }  }^{} \frac{1}{ \|l_{k}\|^{n-1}}\leq \Cl[c]{inductconst}\cdot \|l_{k-1}\|^\delta 
\end{equation}
for some constant $\Cr{inductconst}>0$ and some $0\leq \delta<1$. Indeed, the remaining outer sums then are trivially bounded by summation by parts. 
We split the sum into two parts: one where $ \|l_{k-1}\|\leq \|l_{k}\| \leq \Cr{minima} \|l_{k-1}\|$ and the other $\Cr{minima}\|l_{k-1}\| \leq \|l_{k}\|\leq \Cr{minima} T$. 
Dropping the projection condition on the first sum, it suffices to bound
\begin{equation}
\sum_{ \substack{l_{k} \in \M_{1 \times m}(\mathbb{Z}) \\  \|l_{k-1}\|\leq \|l_{k}\| \leq \Cr{minima} \|l_{k-1}\|  }  }^{} \frac{1}{ \|l_{k}\|^{m}}
+ 
\sum_{ \substack{l_{k} \in \M_{1 \times m}(\mathbb{Z}) \\ \Cr{minima}\|l_{k-1}\| \leq \|l_{k}\| \leq \Cr{minima} T  \\ \|\pi_{k-1}(l_{k})\| \leq \Cr{minima} \|l_{k-1}\| }  }^{} \frac{1}{ \|l_{k}\|^{m}}
\label{eq:two_terms}
\end{equation}
Here, the first sum appears only when $\Cr{minima} > 1$.
Let us focus on the first term. We use~\cref{le:domain_dimension_bound} and 
this gives us that up to some constant, the first term is bounded by 
\begin{equation}
    \log (\Cr{minima} \|l_{k-1}\| )- \log \|l_{k-1}\| .
\end{equation}

For the second sum, we use~\cref{le:domain_dimension_bound} to bound the sum up to a constant by
\begin{equation}
    \|l_{k-1}\| \int_{\|l_{k-1}\|}^{ \Cr{minima} T }  x^{-2} \diff x
\end{equation}
This is because, using now the bound on $\|\pi_{k-1}(l_{k})\|$, we have
\begin{equation}
    \card \{ l_{k
 \in\M_{1 \times m}(\mathbb{Z}) : \|\pi_{k-1}(l_{k})\| \leq \Cr{minima} \|l_{k-1}\| , \|l_{k}\| \leq X
}\} \leq \Cl[c]{finconst} \|l_{k-1}\| X^{m-1},
\end{equation}
for some $\Cr{finconst}>0$. The second sum is therefore bounded. 
Both the terms in~\cref{eq:two_terms} are therefore bounded, up to a constant, by $\log \|l_{k-1}\|$, as desired. The result follows.
\end{proof}

\section{Lifts of codes}
\label{se:lifts_of_codes}

Our goal in this section is to prove the discretized integral formula in \cref{th:higher_moments} using~\cref{th:main}. 

Let $g: \KR^{n} \rightarrow  \mathbb{R}$ be a test function satisfying~\cref{hy:admissible}. 
For any integer $1 \leq s \leq n$ and a prime ideal $\mathcal{P} \subseteq \OK$, denote 
\begin{equation}
    \mathcal{L}(\mathcal{P},s)
    = \{ \tfrac{1}{T_{\mathcal{P}}} \Lambda \mid (\mathcal{P} \mathcal{O}_{K})^{n} \subseteq \Lambda \subseteq \OK^{n} , \Lambda/(\mathcal{P} \OK) ^{n} \in \Gr(s,(\OK/\mathcal{P})^{n}) \},
    \label{eq:def_of_L}
\end{equation}
where after setting
\begin{equation}
    T_{\mathcal{P}} = \N(\mathcal{P})^{\left(1-\frac{s}{n} \right)\frac{1}{d}},
\end{equation}
we get that all the lattices in $\mathcal{L}(\mathcal{P},s)$ to have the same covolume as $\OK^{n} \subseteq \KR^{n}$, which is $1$.
Now we begin considering our object of interest: the average of lattice sum functions. Observe that
\begin{equation}
    \frac{1}{\card \mathcal{L}(\mathcal{P},s)}
    \sum_{ \Lambda \in \mathcal{L}(\mathcal{P},s)} 
    \left(\sum_{v \in \Lambda} g(v)\right)^{m}  = 
    \frac{1}{\card \mathcal{L}(\mathcal{P},s)} 
    \sum_{ \Lambda \in \mathcal{L}(\mathcal{P},s)} 
    \left( \sum_{v \in \Lambda^{m}} f(v)\right),
    \label{eq:nthmoment}
\end{equation}
where $f(v_1,\dots,v_m)= g(v_1)g(v_2)\dots g(v_m)$.
We perform some manipulations on this sum. 
Letting $\ind(P)$ denote the indicator function of a proposition $P$ and writing $x=(x_{1},\dots,x_{m}) \in (\OK^{n})^{m}$, we have that 
\begin{align}
& 
\frac{1}{\card \mathcal{L}(\mathcal{P},s)} \sum_{\Lambda \in \mathcal{L}(\mathcal{P},s)}\left( \sum_{v \in \Lambda^{m}} f(v)\right)  \\
& = 
\sum_{x \in \mathcal{O}_K^{n \times m} }
f(  \tfrac{1}{T_{\mathcal{P}}}  x)
\left(
    \frac{1}{\card \mathcal{L}(\mathcal{P},s)} 
    \sum_{\substack{  S \subseteq k_\mathcal{P}^{n} \\  S \simeq k_{\mathcal{P}}^{s} } } 
    \ind\left(
        \spant
        ( \pi_{\mathcal{P}}(x_1),\dots,\pi_{\mathcal{P}}(x_m)) \subseteq S
\right) \right), 
\label{eq:prepoission}
\end{align}
where $k_{\mathcal{P}} = \OK/\mathcal{P}$ and $\pi_\mathcal{P}: \OK^{n} \rightarrow  k_\mathcal{P}^{n}$ is the ``reduction modulo $\mathcal{P}$'' map. 

The inner sum is just the probability of a random subspace $S \subseteq k_\mathcal{P}^{n}$ of fixed dimension $s$ containing some given set of points $x_1,x_2,\dots,x_m \in k_\mathcal{P}^{n}$. This probability, other than depending on $\mathcal{P},s$, depends only on the $k_\mathcal{P}$-dimension of the subspace generated by $\pi_\mathcal{P}(x_1),\dots,\pi_\mathcal{P}(x_m)$. 
This dimension equals the rank of $\pi_\mathcal{P}(x) \in \M_{ n \times m }(k_\mathcal{P})$ 
which is certainly less than the rank of 
$x\in \M_{n \times m}(\mathcal{O}_K) \subseteq \M_{n \times m }(K)$. So we can split our sum into 
\begin{align}
& = 
\sum_{k=0}^{\min(n,m)}\sum_{\substack{x \in \M_{n \times m }(\mathcal{O}_K) \\ \rank(x) = k}}
\frac{  f( \tfrac{1}{T_{\mathcal{P}}} x)}{\card \mathcal{L}(\mathcal{P},s)} 
\left(
    \sum_{\substack{  S \subseteq k_\mathcal{P}^{n} \\  S \simeq k_{\mathcal{P}}^{s} } } 
    \ind\left(
        \spant
        ( \pi_{\mathcal{P}}(x_1),\dots,\pi_{\mathcal{P}}(x_m)) \subseteq S
\right) \right). 
\label{eq:not_even_the_final_form}
\end{align}

Given $x \in \M_{n \times m }(\mathcal{O}_K)$, we might for some $\mathcal{P}$ encounter a ``rank-drop'' phenomenon, that is $\rank\left(\pi_\mathcal{P}(x)\right) < \rank(x)$. However, the good news is that the matrices $x$ where this rank-drop happens can be ``pushed away'' from the support of $f$, as the following lemma shows.

\begin{lemma}
    Suppose that $x \in M_{n \times m }\left(\mathcal{O}_{K}\right)$ is a matrix with $\rank(x) = k \geq 1$ and $\mathcal{P}$ is a prime ideal in $\mathcal{O}_K$ such that $\rank(\pi_\mathcal{P}(x)) < k$. 
    Then, for any Euclidean norm $\| \cdot \|: M_{n \times m }(K_\mathbb{R})  \rightarrow \mathbb{R}_{\ge 0}$, there exists some $\Cl[c]{push}>0$ 
    depending on $K,\|\cdot\|,n,m$ and independent of $k$ and $\mathcal{P}$ such that $x$ must satisfy
    \begin{equation}
        \|x\| \ge \Cr{push} \N(\mathcal{P})^{\frac{1}{k[K:\mathbb{Q}]}}
    \end{equation}
    \label{le:rankdrop}
\end{lemma}
\begin{proof}
    By choosing a $\mathbb{Z}$-basis of $\OK$, we can embed $\iota:\mathcal{O}_K \hookrightarrow M_{[K:\mathbb{Q}]}(\mathbb{Z})$ as a subring of the square integer matrices of size $[K:\mathbb{Q}]$. 
    Without loss of generality, we assume that the norm $\|\cdot \|$ is the Euclidean norm via the embedding 
    $$\iota: M_{n\times m }(\OK)\hookrightarrow M_{n[K:\mathbb{Q}] \times m[K:\mathbb{Q}]}(\mathbb{Z}) \subseteq \mathbb{R}^{nm[K:\mathbb{Q}]^{2}}.$$ 
    Since $\rank(x)=k$, we know that there exists a non-singular $k \times k$ matrix $a \in \M_{k}(\OK)$ 
    appearing as a submatrix in $x$. 
    However for some such matrix $a$, one must have
    $ 0 \neq \det a \in \mathcal{P}$ otherwise there is no rank-drop modulo $\mathcal{P}$. 
    Therefore, we get that 
    $\N(\mathcal{P}) \mid \N(\det a).$
    Since we know that $0 \neq |\det(\iota(a))| \ge \N(\mathcal{P})$, at least one non-zero integer appearing in the matrix entries of $\iota(a)$ would have absolute value $\geq \Cl[c]{push2} \N(\mathcal{P})^{\frac{1}{k[K:\mathbb{Q}]}}$ for some $\Cr{push2}>0$ independent of $\mathcal{P}$. This produces the same lower bound on the Euclidean norm of $\iota(a)$ up to a constant, and similarly also for $\iota(x)$.
\end{proof}

\begin{lemma}
    \label{le:counting}
    Suppose $y_1, y_2,\dots,y_k \in k_{\mathcal{P}}^{n}$ are linearly independent vectors (over $k_{\mathcal{P}}$). Then the following holds:
    {\small\begin{equation}
            \frac{  1}{\card \mathcal{L}(\mathcal{P},s)} 
            \left(
                \sum_{\substack{  S \subseteq k_\mathcal{P}^{n} \\  S \simeq k_{\mathcal{P}}^{s} } } 
                \ind\left(
                    \spant
                    ( y_1,y_2,\dots,y_k) \subseteq S
            \right) \right) = \begin{cases}
                0 &  \text{ if $s<k$}\\
                \N(\mathcal{P})^{-k(n-s)}\cdot (1+ \varepsilon_{1}) &  \text{ if $s \geq k$.}
            \end{cases}
    \end{equation}}
    where the error term $|\varepsilon_{1}| \leq \Cl[c]{normerror} \N(\mathcal{P})^{-1}$ for some $\Cr{normerror}>0$ not depending on $\mathcal{P}$.
\end{lemma}
\begin{proof}
    The case with $s < k $ is clear. 
    In general for a finite field of size $q$, the number of $u$-dimensional subspaces in a $t$-dimensional vector space is the cardinality of the Grassmannian $\Gr(u,\mathbb{F}_q^t)$ given by a polynomial in $q$. The leading terms of this polynomial are
    \begin{align}
        \frac{(q^t-1)(q^t-q)\cdots (q^t-q^{u-1})}{(q^u-1)(q^u-q)\cdots (q^u-q^{u-1})}
 & = q^{u(t-u)} + \Cl[c]{normerror2} q^{u(t-u)-1} + \dots ,\\
 & =  q^{u(t-u)} ( 1 + \varepsilon_{1}),
    \end{align}
    where $\varepsilon_{1}$ is an error term as given in the statement.
    In our case, $q=\card k_\mathcal{P}= \N(\mathcal{P})$. 

    Up to change of variables, the numerator counts the number of $(s-k)$-dimensional subspaces in a $(n-k)$-dimensional space and therefore has cardinality $\card \Gr(s-k,\mathbb{F}_q^{n-k})$. 
    This is sufficient to get our result.
\end{proof}

\begin{theorem}
    \label{th:higher_moments}
    Take $n \ge 2$, $m \in \{1,\dots,n-1\}$ and choose $s$ as either $n-1$, or any number in $\{m,m+1,\dots,n-1\}$ that satisfies
    \begin{equation}
        1-\frac{s}{n} < \frac{1}{m}.
    \end{equation}

    Let $f:K_{\mathbb{R}}^{n \times m} \rightarrow \mathbb{R}$ be a function satisfying \cref{hy:admissible}. 
    With $\mathcal{L}(\mathcal{P},s)$ defined as in \cref{eq:def_of_L}, 
    we have that as $\N(\mathcal{P}) \rightarrow \infty$
    {\small \begin{equation}
            \frac{1}{\card \mathcal{L}(\mathcal{P},s)} 
            \sum_{ \Lambda \in \mathcal{L}(\mathcal{P},s)} 
            \left( \sum_{v \in \Lambda^{m}} f(v)\right)
            \rightarrow
            \sum_{k=0}^{{m}}\sum_{\substack{D \in \M_{k \times m }(K) \\ 
                    \rank(D) = k \\ 
            D \text{ row reduced echelon}}}
            \mathfrak{D}(D)^{-n} \int_{x \in K_\mathbb{R}^{n \times k }} f(x D ) dx,
            \label{eq:right_side_converge}
    \end{equation}}
    where $\mathfrak{D}(D)$ is 
    as defined in~\cref{de:denonimator}. Here, the term at $k=0$ is just $f(0)$.
\end{theorem}
\begin{proof}
    From the discussion above, we arrive at \cref{eq:not_even_the_final_form}, and it remains to consider
    \begin{align}
        \sum_{k=0}^{m}\sum_{\substack{x \in \M_{n \times m }(\mathcal{O}_K) \\ \rank(x) = k}}
        \frac{  f( \tfrac{1}{T_{\mathcal{P}}} x)}{\card \mathcal{L}(\mathcal{P},s)} 
        \left(
            \sum_{\substack{  S \subseteq k_\mathcal{P}^{n} \\  S \simeq k_{\mathcal{P}}^{s} } } 
            \ind\left(
                \spant
                ( \pi_{\mathcal{P}}(x_1),\dots,\pi_{\mathcal{P}}(x_m)) \subseteq S
        \right) \right).
    \end{align}
    Note that here $T_{\mathcal{P}} =  \N(\mathcal{P})^{\left(1-\frac{s}{n}\right)\frac{1}{d}}$. The rank $k$ ranges within $\{0,1,\dots,m\}$ since $\min(n,m)=m$. Also, since $s \ge m$, we expect the quantity in parentheses to be nonzero from~\cref{le:counting}.\par
    We recall that $\M_{n \times m }(K_{\mathbb{R}})$ has the Euclidean measure given by $n \cdot  m$ copies of the quadratic form coming from (\ref{eq:norm}).
    When $k\geq 1$, we know from Lemma \ref{le:rankdrop} that we will encounter a rank-drop mod $\mathcal{P}$ only if for some predetermined constant $\Cr{push}>0$
    \begin{align}
& \| x \|  \ge \Cr{push} \N(\mathcal{P})^{\frac{1}{kd}}  \\
        \Rightarrow &  \|\tfrac{1}{T_{\mathcal{P}}} x\| \ge \Cr{push} \N(\mathcal{P})^{\frac{1}{d} \cdot \left( \frac{1}{k} - \left(1 - \frac{s}{n}\right)\right)}.
    \end{align}
    Since 
    \begin{equation}
        \frac{1}{k} - \left(1-\frac{s}{n} \right) \ge  \frac{1}{m} - \left(1- \frac{s}{n}\right) > 0,
    \end{equation}
    for a large enough value of $\N(\mathcal{P})$ we have that all the matrices of 
    $x \in \M_{n  \times m }\left(\mathcal{O}_K\right)$ where rank-drop could happen are outside the support of $f$. 
    Let us assume that $\N(\mathcal{P})$ is large enough for this to hold. 
    Hence whenever $f(\tfrac{1}{T_{\mathcal{P}}} x)$ is non-zero, the span of $ \pi_\mathcal{P}(x_1) , \pi_\mathcal{P}(x_2) ,\dots,  \pi_\mathcal{P}(x_m) $ is of the same $k_\mathcal{P}$-dimension as the rank of $x$. 
    Using Lemma \ref{le:counting}, we can rewrite our sum as
    {\small\begin{align}
            \sum_{k=0}^{m}\sum_{\substack{x \in \M_{n \times m }(\mathcal{O}_K) \\ \rank(x) = k}}
            \frac{  f( \tfrac{1}{T_{\mathcal{P}}} x)}{ \N(\mathcal{P})^{k(n-s)} } 
            & =     \sum_{k=0}^{m} \Big(\varepsilon_{\mathcal{P},k}+T_{\mathcal{P}}^{-knd}
                        \sum_{\substack{x \in \M_{n \times m }(\mathcal{O}_K) \\ \rank(x) = k}} 
                        {  f( \tfrac{1}{T_{\mathcal{P}}} x)} \Big),\\
            & =  \varepsilon_{2} +   \sum_{k=0}^{m} \Big(T_{\mathcal{P}}^{-knd}
                        \sum_{\substack{x \in \M_{n \times m }(\mathcal{O}_K) \\ \rank(x) = k}} 
                        {  f( \tfrac{1}{T_{\mathcal{P}}} x)} \Big),\\
    \end{align}}
    where $|\varepsilon_{2}| \leq \Cl[c]{grassmannian} \N(\mathcal{P})^{-1}$ for some $\Cr{grassmannian}>0$ that does not depend on $\mathcal{P}$.

    The result follows as $T_{\mathcal{P}} \rightarrow \infty$ due to $\N(\mathcal{P}) \rightarrow  \infty$ using \cref{th:main} and ~\cref{eq:value_of_c}.
\end{proof}

\begin{remark}
    %
    It may be possible that equidistribution results for Hecke points as in \cite{EOL01} could imply the equidistribution of $\mathcal{L}(\mathcal{P},s)$ in the relevant moduli space of $\OK$-modules. 
    Then, as $\N(\mathcal{P}) \rightarrow \infty$, one could gain via \cref{th:higher_moments} yet another proof of the Rogers integral formula studied in 
    \cite{K19,  GSV23}. 
    However, to finish such an argument, we would need equidistribution of these points in a pointwise sense to handle non-compactly supported functions, which seems to be either missing from the literature or out of reach.
\end{remark}
\begin{remark}
One can find the rate of convergence from \cref{th:main}.
\end{remark}

\bibliographystyle{amsalpha}
\bibliography{authfile} 
\end{document}